\documentclass[11pt]{amsart}

\usepackage{amsfonts,amssymb,amsmath,graphicx}
\usepackage[mathscr]{eucal}

\newtheorem{thm}{Theorem}
\numberwithin{thm}{subsection}
\numberwithin{equation}{subsection}

\newtheorem{scholium}[thm]{Scholium}
\newtheorem{propty}[thm]{}

\newtheorem{prop}[thm]{Proposition}
\newtheorem{cor}[thm]{Corollary}
\newtheorem{lem}[thm]{Lemma}
\newtheorem{rem}[thm]{Remark}
\newtheorem{remarks}[thm]{Remarks}
\newtheorem{hist}[thm]{Historical note}

   \def\Int{\operatorname{Int}} 
   \def\Im{\operatorname{Im}} 
   \def\Re{\operatorname{Re}} 
   \def\Aut{\operatorname{Aut}} 
   \def\card{\operatorname{card}} 
   \def\pr{\operatorname{pr}} 
   \def\top{\operatorname{top}} 
   \def\bot{\operatorname{bot}} 
   \def\a{\alpha} \def\b{\beta} 
   \def\e{\varepsilon} \def\g{\gamma} \def\k{\kappa} 
   \def\phi{\varphi}  \def\s{\sigma} \def\t{\tau} 
\def\twobars#1#2#3#4{\vcenter{\hrule height.#1pt width#2pt
                               \vskip#3pt
                               \hrule height.#1pt width#2pt
                               \vskip#4pt}}
\def\stroke#1#2#3{\vrule height#1pt width.#2pt depth#3pt}
\def\bdconnsum{\hskip2pt
             \stroke84{-1.55}\twobars5331\stroke{5.5}41
             \hskip2pt}
     \def\Bd{\partial} 
     \def\arep{{\vec a}} 
         \def\brep{{\vec b}} 
         \def\crep{{\vec c}} 
         \def\drep{{\vec d}} 
     \def\concat{\smallfrown} 
     \def\eedge{\mathbf e}\def\fedge{\mathbf f}\def\gedge{\mathbf g}
     \def\iso{\cong} 
     \def\numk{\mathbf k}
     \def\numn{\mathbf n}
     \def\C{\mathbb C} \def\N{\mathbb N} 
     \def\R{\mathbb R} 
     \def\Akn{A(K,n)}  
     \def\AKN#1{A(K_#1,n_#1)} \def\AON#1{A(O_#1,n_#1)} 
     \def\Col#1#2{\operatorname{col}_{\,#2}{(#1)}}
     \def\h#1#2{h^{\scriptscriptstyle(#1)}_{#2}} 
     \def\Nb#1#2{N_#1(#2)} 
     \def\sub{\subset} \def\emptyset{\varnothing}
     \def\Vert{\mathsf V} 
     \def\Edge{\mathsf E} 
     \def\Endpt{\mathsf V_{\mathrm{end}}} 
     \def\Xend{X_{\mathrm{end}}} 
     \def\star{\ast}
     \def\plumb#1{\ast_#1} 
     \def\topplumb#1{\thinspace{\overline\ast\thinspace}_#1 \thinspace}
     \def\botplumb#1{\thinspace{\underline\ast\thinspace}_#1 \thinspace}
     \def\PLUMB{\mathop{\star}\nolimits}
     \def\TOPPLUMB{\mathop{\overline\star\thinspace}\nolimits}
     \def\BOTPLUMB{\mathop{\underline\star\thinspace}\nolimits}
     \def\As{\mathscr A}  
       
     \def\Is{\mathscr I}

     \def\Ss{\mathscr S} \def\Ts{\mathscr T} 
       
     \def\Ys{\mathscr Y}  
     \def\LHP{\C_{-}} 
     \def\UHP{\C_{+}} 
     \def\commute#1#2{\lbrack#1,#2\rbrack} 
     \def\YBax#1#2
         {\medspace\llless\negmedspace#1,#2\negmedspace\gggtr\medspace}
\begin{document}

\title[Hopf plumbing and homogeneous braids]
{Hopf plumbing, arborescent {S}eifert surfaces, 
\hbox{baskets,} espaliers, and homogeneous braids}
\author{Lee Rudolph}
\address{Department of Mathematics, Clark University, 
Worcester MA 01610 USA}  
\thanks{Partially supported by UNAM, CAICYT, NSF (DMS-8801915, 
DMS-9504832), and CNRS}
\email{lrudolph@black.clarku.edu}            

\begin{abstract}
Four constructions of Seifert surfaces---Hopf and arborescent 
plumbing, 
basketry, 
and \hbox{$\Ts$-bandword} handle decomposition---are described, 
and some interrelationships expounded, e.g.: arborescent Seifert 
surfaces are baskets; Hopf-plumbed baskets are precisely 
homogeneous $\Ts$-bandword surfaces. 
\end{abstract}

\maketitle

\section{Introduction; statement of results}\label{intro}

Let $A(O,n)\sub S^3$ denote an $n$-twisted unknotted annulus.
A Seifert surface $S$ is \emph{Hopf-plumbed} (see \S\ref{plumbing})
if $S=D^2$ or if $S=S_0\plumb{\a} A(O,\mp 1)$ can be constructed 
by plumbing a positive or negative Hopf annulus $A(O,\mp 1)$ 
to a Hopf-plumbed surface $S_0$ along a proper arc $\a\sub S_0$; 
a Hopf-plumbed surface is a fiber surface, 
\cite{Stallings,Harer,Morton}.
A Seifert surface $S$ is a \emph{basket} 
(see \S\ref{basketry}) 
if $S=D^2$ or 
if $S=S_0\plumb\a A(O,n)$ can be constructed by plumbing $A(O,n)$ 
to a basket $S_0$ along a proper arc $\a\sub D^2\sub S_0$, 
\cite{Rudolph98,Boileau-Rudolph}; a fundamental theorem 
of Gabai \cite{Gabai83b,Gabai85} implies that a basket is a fiber 
surface iff\/ it is Hopf-plumbed.
A Seifert surface $S$ is \emph{arborescent} 
(see \S\ref{arborescence}) if 
$S=D^2$, or if $S=A(O,n)$, or if $S=S_0\plumb\a A(O,n)$ can be
constructed by plumbing $A(O,n)$ to an arborescent Seifert 
surface $S_0$ along a transverse arc $\a$ of an annulus
plumband of $S_0$, \cite{Conway,Siebenmann,Gabai86}; 
an arborescent Seifert surface is a basket 
(Proposition \ref{arborescent surface is basket}).

In \S\ref{espaliers}, I define the \emph{$\Ts$-generators} of a braid 
group corresponding to a tree $\Ts\sub\C$.  
Call $\Ts$ an \emph{espalier} if each edge $\eedge$ of $\Ts$
is properly embedded in $\LHP:=\{z\in\C:\Im{z}\le 0\}$
and $\Re|\,\eedge:\eedge\to\R$ is injective.
For an espalier $\Ts$, 
words $\brep$ in the \hbox{$\Ts$-generators} correspond nicely
to \hbox{\emph{$\Ts$-bandword surfaces}} $S(\brep)\sub\R^3\sub S^3$.
(If $\Is_\numn$ is an espalier with 
an edge from $p$ to $p+1$, $1\le p<n$,
then the $\Is_\numn$-generators of $B_n$ are the standard generators 
$\s_1,\dots,\s_{n-1}$, and $S(\brep)$ is the output of Seifert's 
algorithm \cite{Seifert} applied to the closed braid diagram of $\brep$.)  
If $\brep$ is \emph{homogeneous}, then $S(\brep)$ 
is a fiber surface, \cite{Stallings,Rudolph83b,Rudolph88}; 
in particular $S(\brep)$ is incompressible and connected. 
In \S\ref{characterization}, I generalize (and 
make more precise) the fact that 
a homogeneous $\Is_\numn$-bandword surface is Hopf-plumbed,
and prove a converse.  

\noindent
\textbf{Main Theorem.}\emph{
\emph{(1)} If $S$ is any Hopf-plumbed basket (for instance, an 
arborescent fiber surface), then there is an espalier
$\Ts$ and a homogeneous $\Ts$-bandword $\brep$ such that 
$S$ is isotopic to $S(\brep)$.
\emph{(2)} If $\Ts$ is any espalier, and 
$\brep$ is a $\Ts$-bandword,
then the following are equivalent:
\emph{(A)} $\brep$ is homogeneous;
\emph{(B)} $S(\brep)$ is a Hopf-plumbed basket;
\emph{(C)} $S(\brep)$ is a fiber surface;
\emph{(D)} $S(\brep)$ is incompressible and connected.  
}%

General notations, definitions, and conventions are established in
\S\ref{prelims}.  A number of more or less well-known results 
about plumbing are
collected in \S\ref{plumbing}; the exposition and notation 
there are adapted to the requirements of 
\S\S\ref{basketry and arborescence}-\ref{characterization}.

\section{Preliminaries}\label{prelims}
\subsection{Miscellany}\label{misc}
Both $A:=B$ and $B=:A$ define $A$ as meaning $B$.
The symbol $\square$ signals either the 
end or the omission of a proof.
Projection on the $i$th factor of a cartesian product 
is denoted by $\pr_i$.  
For $n\in\N:=\{0,1,2,\dots\}$, $\numn:=\{1,\dots,n\}$.
Let $s,t,s',t'\in \R$, $s\ne t$, $s'\ne t'$.
Say $\{s,t\}$ and $\{s',t'\}$ \emph{touch}
if $\card(\{s,t\}\cap\{s',t'\})=1$, and 
\emph{link} (resp., \emph{unlink})
iff\/ they do not touch and the cross-ratio
$((s-s')(t-t'))/((s-t')(s'-t))$
is positive (resp., negative).

\subsection{Spaces}\label{spaces}
Spaces, maps, etc., are piecewise smooth.
Isotopies are ambient unless otherwise noted.
That $X$ is isotopic to $Y$ is denoted $X\iso Y$,
or by $X\iso_Z\! Y$ if there is such an isotopy fixing $Z$ pointwise.
Manifolds may have boundary and are oriented unless otherwise noted;
in particular, $\R$, $\C^n$, and 
$S^{2n-1} \sub \C^n$
have standard orientations, as does $\R^3$ which is identified 
with the complement of a point $\infty\in S^3$.  
Identify $\C=\R+\sqrt{-1}\,\R\supset\R$ with $\R^2$
and write $\Re$ (resp., $\Im$) for $\pr_1$ (resp., $\pr_2$).  
Then also $\C\times\R$ is identified with $\R^3$.
Write $\LHP:=\{w\in\C:\Im{w}\le 0\}$ (resp., 
$\UHP:=\{w\in\C:\Im{w}\ge 0\}$) for the closed lower 
(resp., upper) half-plane.

If $M$ is a manifold, 
then $-M$ denotes $M$ with 
its orientation reversed (and, where notation
requires it, $+M$ denotes $M$).  
For a suitable subset $Q\sub M$, $\Nb MQ$ denotes a closed 
regular neighborhood of $Q$ in $(M,\Bd M)$.
For a suitable codimension-$1$ submanifold $Q\sub M$ (resp., 
submanifold pair $(Q,\Bd Q)\sub(M,\Bd M)$), a \emph{collaring} 
is an orientation-preserving embedding 
\hbox{$Q\times[0,1]\to M$} 
(resp., $(Q,\Bd Q)\times[0,1]\to (M,\Bd M)$) 
extending $1_Q=1_{Q\times\{0\}}$; a \emph{collar} 
of $Q$ in $M$ (resp., of $(Q,\Bd Q)$ in $(M,\Bd M)$) is 
the image $\Col{Q}{M}$ (resp., $\Col{Q,\Bd Q}{(M,\Bd M)}$) 
of a collaring.
The \emph{push-off} of $Q$ determined by a collaring of $Q$ or $(Q,\Bd Q)$,
denoted by $Q^{+}$, is the image by the collaring of $Q\times\{1\}$ with
the orientation of $Q$; let $Q^{-}:=-Q^{+}$ (so that $Q$ and $Q^{-}$ are 
oriented submanifolds of the boundary of $\Col{Q}{M}$).

An \emph{arc} is a manifold diffeomorphic to $[0,1]$.
An \emph{edge} is an unoriented arc.
A \emph{surface} is a compact $2$-manifold no component of which 
has empty boundary.
An arc or edge $\a$ in a $2$-manifold $S$ is \emph{proper} 
(resp., \emph{boundary}) if $\Bd\a=\a\cap\Bd S$ (resp., $\a\sub\Bd S$).
Given ordered index sets $X, U$, 
and an ordered
handle decomposition
\begin{equation}\label{handle decomposition}
S=\bigcup_{x\in X}\h0x\cup \bigcup_{u\in U}\h1u 
\end{equation}
of a surface $S$, write $S^{(0)}:=\bigcup_{x\in X}\h0x$,
$S^{(1)}:=\bigcup_{u\in U}\h1u$; it is understood
that if $u<v$ then the attaching region
$\h1{u}\cap\Bd S^{(0)}$ is disjoint from $\h1{v}\cap\Bd S^{(0)}$.
A \emph{core} (resp., \emph{transverse}) arc of a $1$-handle 
$\h1{}$ is any proper arc $\k(\h1{})$ (resp., $\t(\h1{})$)
which joins interior points of the two components 
of the attaching region of $\h1{}$ (resp., the complement
in $\Bd \h1{}$ of the attaching region of $\h1{}$).

\subsection{Seifert surfaces}\label{Seifert surfaces}
A \emph{Seifert surface} is a surface $S\sub S^3$.  
Let $\top(S):=\Col{S}{S^3}$, $\bot(S):=\top(-S)$.
A \emph{link} $L$ is the boundary of a Seifert 
surface.  A \emph{knot} is a connected link; a knot
$O$ which is the boundary of a disk $D^2\sub S^3$ is an \emph{unknot}.
If $K$ is a knot, then $\Akn$ denotes any Seifert surface $A$
(necessarily an annulus) such that $K\sub\Bd A$,
$A$ is a collar of $K$ in $A$, and the linking number
in $S^3$ of $K$ and $K^{+}$ is $n$.  Since clearly $K^{-}\iso-K$, 
so also $\Akn\iso A(-K,n)$; further, $-\Akn\iso \Akn$.  
A \emph{transverse arc} of
$\Akn$ is any proper arc $\t(\Akn)\sub\Akn$ from $K$ to $K^{-}$;
$\t(\Akn)$ is unique up to isotopy on~$\Akn$.

\subsection{Incompressible surfaces; fiber 
surfaces}\label{incompressible surfaces}
Let $S$ be a Seifert surface.  Say that $D^2\sub S^3$ is a 
\emph{top-compression disk} for $S$, 
and call $S$ \emph{top-compressible}, if 
$\Bd D^2 =  D^2\cap S$, 
$\Bd D^2$ bounds no disk on $S$,
and $D^2\cap\top(S)=\Col{\Bd D^2}{D^2}$; 
a top-compression disk for $S$ can always be taken 
to be disjoint from $\Int\bot(S)$.
It is a well known consequence of the Loop Theorem
that $S$ is top-compressible 
iff $\pi_1(S_0)\to\pi_1(S^3\setminus S)$ 
(induced by the inclusion $S^{+}\hookrightarrow S^3\setminus S$
of a push-off) is not injective for some component $S_0$ of $S$.
Call $S$ \emph{compressible} if either $S$ or $-S$
is top-compressible, and \emph{incompressible} otherwise.  
(It can happen that exactly one of $S$, $-S$ is top-compressible,
cf.\ \cite{Gabai86}.)
Call $S$ a \emph{fiber surface} 
(and $\Bd S$ a \emph{fibered link}) if there 
is a fibration $\phi:S^3\setminus \Bd S\to S^1$ 
such that $\Int S=\phi^{-1}(1)$ and for all $\zeta\in S^1$,
the closure of $\phi^{-1}(\zeta)$ is a surface with 
boundary $\Bd S$.  
By \cite{Stallings} and \cite{Gabai83},
$S$ is a fiber surface iff $S$ is connected and 
$\pi_1(S)\to\pi_1(S^3\setminus S)$ is bijective.  
Consequently, a fiber surface is incompressible,
and $\Akn$ is a fiber surface iff $K=O$, $n=\mp1$.

\subsection{Trees}\label{trees}
A \emph{tree} is a finite, connected, 
acyclic $1$-complex $\Ts$. 
Let $\Vert(\Ts)$ (resp., $\Endpt(\Ts)$; $\Edge(\Ts)$) 
denote the set of $0$-cells (resp., endpoints; $1$-cells) of $\Ts$.

\section{Plumbing}\label{plumbing} 
This section is a systematic exposition of results about 
those Seifert surfaces that can be constructed from some base 
Seifert surface (not necessarily a disk) by iterated plumbing 
of annuli.  Although most (if not all) of the results are 
more or less well-known, there are a few I have not been able
to find in the literature.  

\subsection{Top- and bottom-plumbing}\label{top and bottom}
\begin{figure}
\begin{center}
\includegraphics{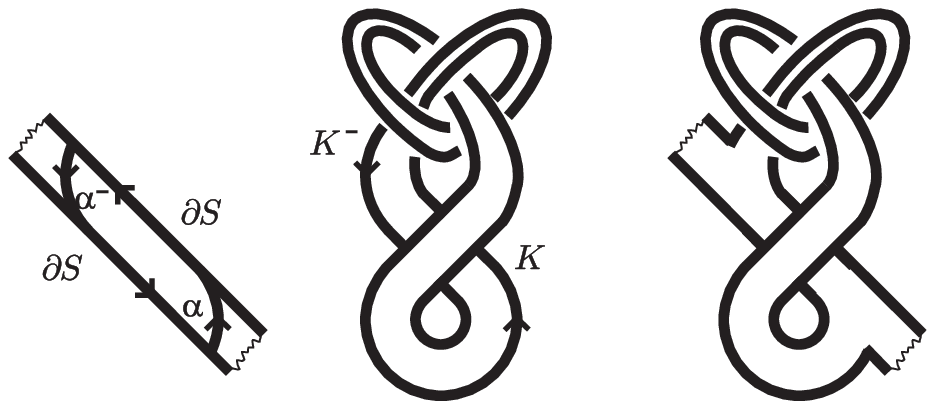}
\end{center}
\caption{A neighborhood of $C_\a$ on $S$; 
an annulus $\Akn$; the top-plumbed Seifert
surface $S\topplumb\a\Akn$.}
\label{Figure 1}
\end{figure}

\begin{figure}
\begin{center}
\includegraphics{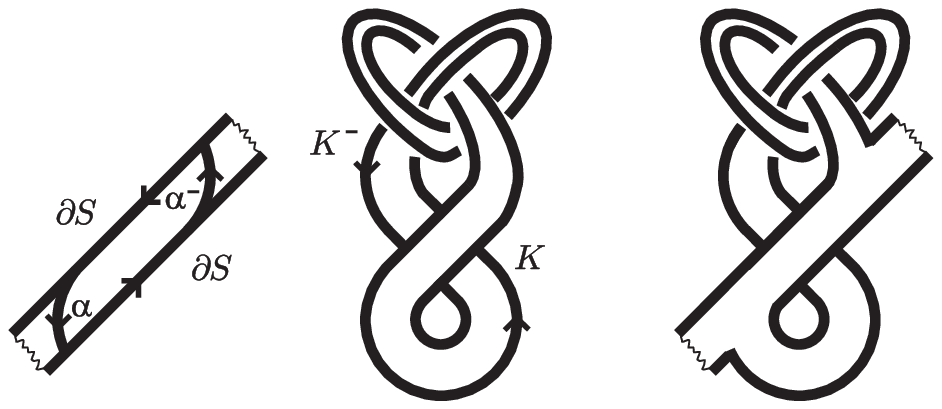}
\end{center}
\caption{Bottom-plumbing $\Akn$ to $S$ along $\a$.}
\label{Figure 2}
\end{figure}

Let $\a\sub S$ be a proper arc on a Seifert surface.
Let $C_\a:=\Col{\a,\Bd\a}{(S,\Bd S)}$.  Then $C_\a$ is a 
\emph{$2$-patch} in the sense of \cite{Rudolph98}, that is, 
a $2$-cell naturally endowed with the structure of a $4$-gon
such that $\Bd C_\a$ is the union of two proper arcs in $S$
(namely, $\a$ and $\a^{-}$) and two boundary arcs in $S$,
say $\g_\a$ and $\g_\a^{-}$, in the cyclic order $\a$, $\g_\a$, 
$\a^{-}$, $\g_\a^{-}$.  Let $D_\a:=\Col{C_\a}{\top(S)}$
(so $D_\a$ is a $3$-cell ``on top'' of $S$, that is, the
positive normal to $S$ along $C_\a = S\cap D_\a\sub\Bd D_\a$ 
points into $D_\a$).  Let $\Akn\sub D_\a$ be an annulus such that
$\Akn\cap\Bd D_\a=C_\a$, $\a\sub K$, $\a^{-}\sub K^{-}$.  Then
$\g_\a\sub\Akn$ is a proper arc, 
$C_\a=\Col{\g_\a,\Bd\g_\a}{(\Akn,\Bd\Akn)}$ is a
$2$-patch on $\Akn$, and 
the union
$S\cup\Akn=:S\topplumb\a\Akn$ 
is a Seifert surface.  
Say $S\topplumb\a\Akn$ is
constructed by \emph{top-plumbing $\Akn$ to $S$ along $\a$}.
Figure~\ref{Figure 1} illustrates this construction;
Figure~\ref{Figure 2} illustrates 
\emph{bottom-plumbing $\Akn$ to $S$ along $\a$},
that is, the construction of 
$S\botplumb\a\Akn := -(-S\topplumb {{{-}\a}} A(-K,n))$ 
from $\a\sub S$ and $\Akn$.
In any case where 
the distinction between $S\topplumb\a\Akn$ and $S\botplumb\a\Akn$ 
can safely be suppressed 
(e.g., if they are isotopic; see \ref{coincidences}--\ref{commuting}), 
each may be denoted 
by $S\plumb\a\Akn$ and said simply to have been 
constructed by \emph{plumbing $\Akn$ to $S$ along $\a$}. 
Given $S':=S\plumb\a\Akn$, call any transverse arc of $\Akn$ which 
is disjoint from $S$ a \emph{transverse arc of $\Akn$ in $S'$},
and denote it by $\t(\Akn\sub S')$.


It is useful to have some notations for iterated plumbing.
Given proper arcs $\a_1\sub S_0$, $\a_2\sub S_1:=S_0\plumb{{\a_1}}\AKN1$, 
\dots, $\a_\mu\sub S_{\mu-1}$, call
\begin{multline}\label{iterated plumbing}
S_\mu=(\dots(S_0\plumb{{\a_1}}{\AKN1}\plumb{{\a_2}}{\AKN2}\dots)%
\plumb{{\a_\mu}}{\AKN\mu}\\
=:S_0\PLUMB_{(\a_1,\dots,\a_\mu)} (\AKN1,\dots,\AKN\mu)
\end{multline}
an \emph{annulus $S_0$-presentation of $S_\mu$} with 
\emph{plumbands} $\AKN{s}$ and \emph{plumbing arcs} $\a_s$;  
further abbreviate the righthand side of (\ref{iterated plumbing}) 
by $S_0\PLUMB_{\vec\a} A(\vec K, \vec n)$ 
(where $\vec\a:=(\a_1,\dots,\a_\mu)$, etc.).  
Call $S$ \emph{annulus $S_0$-plumbed} if it has an 
annulus $S_0$-presentation.

\begin{rem}\emph{ 
The notation (\ref{iterated plumbing})
does not indicate the type (top or bottom) of each plumbing.  
In cases where all plumbings are of one type, more specific notations 
$S_0\TOPPLUMB_{(\a_1,\dots,\a_\mu)} (\AKN1,\dots,\AKN\mu)$
(or $S_0\TOPPLUMB_{\vec \a} A(\vec K, \vec n)$) and 
$S_0\BOTPLUMB_{(\a_1,\dots,\a_\mu)} (\AKN1,\dots,\AKN\mu)$
(or $S_0\BOTPLUMB_{\vec \a} A(\vec K, \vec n)$) 
may be used.
}\end{rem}

\begin{hist}\label{history of plumbing}\emph{
Top-plumbing is a special 
case of \emph{Stallings plumbing} \cite{Stallings} 
(equivalently, \emph{Murasugi sum} \cite{Gabai83b}), 
in the general case of which Seifert surfaces $S_1$, $S_2$ (possibly 
neither an annulus) are attached along a $k_1$-patch and a $k_2$-patch
(possibly with $k_i>2$).  The special case of annulus $D^2$-plumbing
in which all plumbands are unknotted is also 
a special case of \emph{arborescent plumbing} 
\cite{Conway,Siebenmann,Gabai86} (see \S\ref{arborescence}), 
which in general also admits
unknotted M\"obius bands as plumbands (of course 
``top'' has no global meaning where 
unoriented, possibly nonorientable, surfaces are involved).
}\end{hist}

\subsection{Coincidences among isotopy classes of 
plumbed surfaces}\label{coincidences}
Up to isotopy in $S^3$, $S\topplumb\a\Akn$ depends 
only on $S$, $\a$ (up to isotopy on $S$), $K$, $n$,
and the orientations of $S$, $\a$, and $K$.  
In all cases $S\topplumb\a\Akn=S\topplumb{{-\a}}{A(-K,n)}$.
In many cases, no two of the Seifert surfaces
$S\topplumb\a\Akn$, $S\botplumb\a\Akn$, 
$S\topplumb{{{-}\a}}\Akn$, and $S\botplumb{{{-}\a}}\Akn$ 
are isotopic, though each is the union of oriented submanifolds
$S$ and $\Akn$ intersecting in $C_\a$:
if the knot $K$ is not reversible (i.e., $K\not\iso-K$), then 
$S\topplumb\a\Akn$ and $S\topplumb\a A(-K,n)$ may not be 
isotopic; and if neither $S$ nor $\Akn\ne A(O,\mp 1)$ is a fiber surface,
then $S\topplumb\a\Akn$ and $S\botplumb\a\Akn$ may not be isotopic.  
However, various systematic coincidences are worth noting.
Let $S$ be a Seifert surface, $\a\sub S$ a proper arc.

\begin{figure}
\begin{center}
\includegraphics{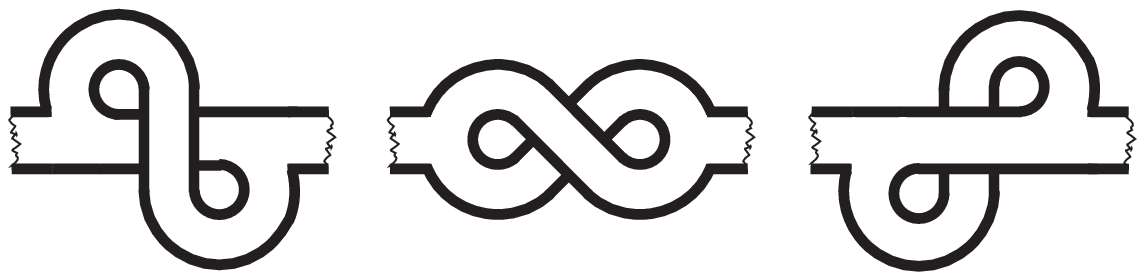}
\end{center}
\caption{Indifference of Hopf annuli to ``top'' and ``bottom'':
an isotopy from
$S\topplumb \a A(O,-1)$ to $S\botplumb \a A(O,-1)$.
}
\label{Figure 3}
\end{figure}

\begin{lem}\label{indifference}
$S\topplumb\a A(O,\mp 1)\iso S\botplumb\a A(O,\mp 1)$.
\end{lem}
\begin{proof} This is a well-known consequence of general facts
about geometric monodromies of fibered links
(see, e.g., \cite{Morton}).  It is also easily 
seen directly (Figure~\ref{Figure 3}). 
\end{proof}

\begin{scholium} If $S$ is a fiber surface, 
then $S\topplumb\a\Akn\iso S\botplumb\a\Akn$.
\end{scholium}

\begin{proof} This follows from the same generalities as 
\ref{indifference}.
\end{proof}

\begin{figure}
\begin{center}
\includegraphics{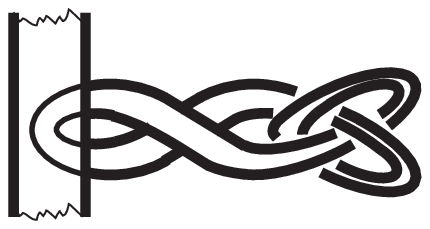}
\end{center}
\caption{$S \bdconnsum \Akn=S\plumb\a\Akn$.}
\label{Figure 4}
\end{figure}

\begin{lem} If $\a$ is boundary-compressible 
(i.e., there is $2$-disk 
$D\sub S$ with $\a\sub\Bd D$, $\Bd D\setminus\Int\a\sub\Bd S$), 
then for any $\Akn$, $S\topplumb\a\Akn\iso S\botplumb\a\Akn$.
\end{lem}
\begin{proof} Each is a boundary-connected 
sum $S \bdconnsum \Akn$ along $\Bd D\setminus\Int\a$ 
and an arc of $K\sub\Bd\Akn$ (Figure~\ref{Figure 4}). \end{proof}

\begin{lem}\label{summing is plumbing}
Any Murasugi sum of 
$S$ with an annulus $\Akn$ 
(not necessarily along $2$-patches)
is isotopic to $S\plumb\a\Akn$ for an appropriate $\a$.
\qed\end{lem}
\begin{figure}
\begin{center}
\includegraphics{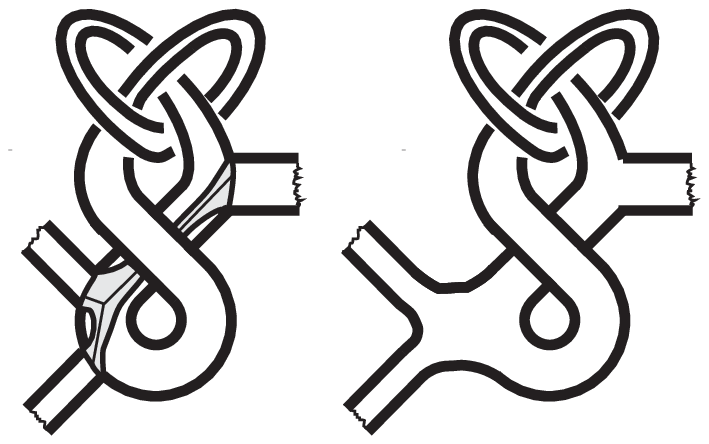}
\end{center}
\caption{A Murasugi sum $S\ast\Akn$ along a $3$-patch,
and an isotopic top-plumbing $S\topplumb\a\Akn$.}
\label{Figure 5}
\end{figure}
\begin{proof} This is a consequence of the paucity of distinct
isotopy classes of $k$-patches on an annulus.  
The example in Figure~\ref{Figure 5} adequately suggests the 
general proof.\end{proof}
\begin{lem} Let $\g\sub S$ be a proper arc 
such that $S=S'\bdconnsum_\g S''$. 
If ${\overline S}_0:=S\topplumb\g A(K_0,n_0)$, 
${{\underline S}\,}_0:=S\botplumb\g A(K_0,n_0)$,
$S_\mu:=S\PLUMB_{\vec\a} A(\vec K, \vec n)$,
then ${\overline S}_0\PLUMB_{\vec\a} A(\vec K, \vec n)$ and
${{\underline S}\,}_0\PLUMB_{\vec\a} A(\vec K, \vec n)$ are well-defined,
and ${\overline S}_0\PLUMB_{\vec\a} A(\vec K, \vec n) \iso_{S_\mu}
{{\underline S}\,}_0\PLUMB_{\vec\a} A(\vec K, \vec n)$.
\end{lem}
\begin{proof} All the plumbing arcs $\a_s$ miss 
$A(K_0,n_0)\setminus S$, so 
${\overline S}_0\PLUMB_{\vec\a} A(\vec K, \vec n)$ and
${{\underline S}\,}_0\PLUMB_{\vec\a} A(\vec K, \vec n)$ are 
well-defined.
Figure~\ref{Figure 6} illustrates the required isotopy.\end{proof}
\begin{figure}
\begin{center}
\includegraphics{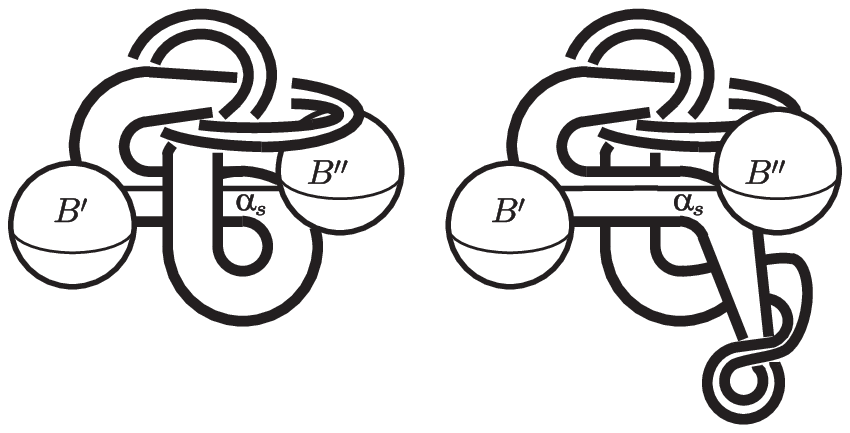}
\end{center}
\caption{Disjoint $3$-balls $B'$ and $B''$ enclose $S'$ and $S''$; %
there is no obstruction to pulling $A(K_0,n_0)$ around from front to back.}
\label{Figure 6}
\end{figure}

\begin{figure}
\begin{center}
\includegraphics{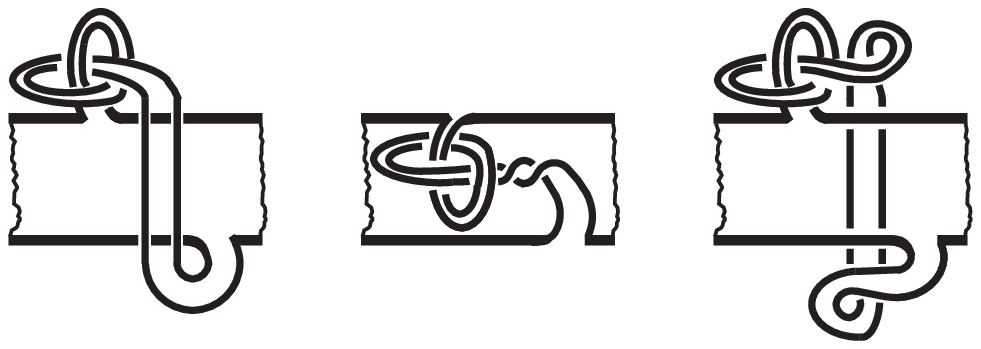}
\end{center}
\caption{An isotopy from $\Bd (S\topplumb\a\Akn)$ 
to $\Bd (S\botplumb\a\Akn)$, pushing 
$\Bd(\Akn\setminus S)$ through $\Akn\cap S$.}
\label{Figure 7}
\end{figure}

\begin{cor}\label{back-to-front}
 For any $\gamma$, $\vec\a$, $A(K_0,n_0)$, 
and $A(\vec\a,\vec n)$,
\begin{displaymath}
\quad\quad\quad
(D^2\topplumb\gamma A(K_\gamma,n_\gamma))\PLUMB_{\vec\a} %
\!A(\vec K, \vec n)\!\iso_{D^2}\!(D^2\botplumb\gamma\Akn)\PLUMB_{\vec\a} %
\!A(\vec K, \vec n).\quad\quad\thickspace\square
\end{displaymath}
\end{cor}
\begin{scholium} There exists a non-ambient isotopy 
in $(D^4,S^3)$ of the pairs 
$(S\topplumb\a\Akn,\Bd (S\topplumb\a\Akn))$ and 
$(S\botplumb\a\Akn,\Bd (S\botplumb\a\Akn))$;
in particular, $\Bd (S\topplumb\a\Akn)\iso \Bd (S\botplumb\a\Akn)$.
\end{scholium}
\begin{proof} Figure~\ref{Figure 7}.
 \end{proof}
\subsection{Commuting plumbings}\label{commuting}
Let $S$ be a Seifert surface, $\a_1, \a_2\sub S$ proper arcs.
\begin{lem}\label{top/bottom} For all $K_1, K_2, n_1$, and $n_2$,
\begin{displaymath}
(S\topplumb{{\a_1}}A(K_1,n_1))\botplumb{{\a_2}}A(K_2,n_2)=
(S\botplumb{{\a_2}}A(K_2,n_2))\topplumb{{\a_1}}A(K_1,n_1).
\end{displaymath}
\end{lem}
\begin{proof} Clear, since (if minimal care is taken)
$\Int\top(S)\cap\Int\bot(S)=\emptyset$.
\end{proof}

\begin{lem} If $|n_1|=1=|n_2|$, then 
$$
S\PLUMB_{(\a_1,\a_2)} (A(O,n_1),A(O,n_2))\iso
S\PLUMB_{(\a_2,\a_1)} (A(O,n_2),A(O,n_1)).
$$
\end{lem}
\begin{proof} Immediate from \ref{indifference} and 
\ref{top/bottom}.
\end{proof}

\begin{figure}
\begin{center}
\includegraphics{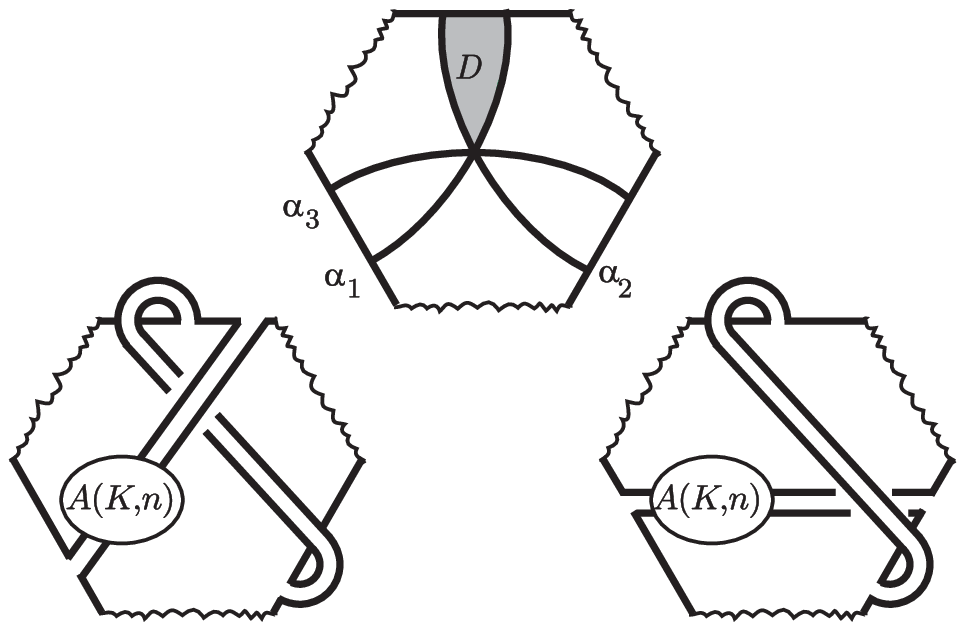}
\end{center}
\caption{One attaching region
of the $1$-handle $\Akn\setminus\Int S$ slides off $\Bd S$,
over the $1$-handle $A(O,1)\setminus\Int S$, and back
to $\Bd S$.}
\label{Figure 8}
\end{figure}

\begin{lem} 
If $\a_1\cap\a_2=\emptyset$, 
then for all $K_1$, $K_2$, $n_1$, and $n_2$,
$$
S\PLUMB_{(\a_1,\a_2)} (A(K_1,n_1),A(K_2,n_2))
\iso
S\PLUMB_{(\a_2,\a_1)} (A(K_2,n_2),A(K_1,n_1)).
$$
\end{lem}
\begin{proof} If (without loss of generality)
$C_{\a_1}$ and $C_{\a_2}$ are taken to be disjoint,
then in fact the two surfaces can be taken to be identical.
\end{proof}
\begin{scholium} 
If $\a_1, \a_2\sub S$ intersect transversely
in a single point $P$, and $D\sub S$ is a $3$-gon 
such that $\Bd D\cap \Bd S$ is an arc and 
$\Bd D\cap \a_s$ is an arc with endpoint $P$,
then there is a proper arc $\a_3\sub S$ 
intersecting each of $\a_1, \a_2$ transversely 
in the single point $P$, such that, for any $\Akn$,
\begin{displaymath}
S\PLUMB_{(\a_1,\a_2)} (\Akn,A(O,1))
=S\PLUMB_{(\a_2,\a_3)} (A(O,1),\Akn),
\end{displaymath}
\begin{displaymath}
S\PLUMB_{(\a_2,\a_1)} (A(O,-1),\Akn)
=S\PLUMB_{(\a_3,\a_2)} (\Akn,A(O,-1)).
\end{displaymath}
\end{scholium}
\begin{proof}
The first case is illustrated in Figure~\ref{Figure 8}; 
the second case is similar.\end{proof}

\subsection{Annulus-plumbed and Hopf-plumbed Seifert 
surfaces}\label{annulus-plumbed}
An annulus $D^2$-plumbed Seifert surface will be called simply 
\emph{annulus-plumbed}.  By \ref{summing is plumbing}, 
any Seifert surface which can 
be constructed from $D^2$ by iterated Murasugi sums with annuli 
is in fact annulus-plumbed.
A \emph{Hopf $S_0$-presentation} is an annulus $S_0$-presentation 
(\ref{iterated plumbing}) such that each plumband is a Hopf 
annulus $A(O,\mp1)$.  A Hopf $D^2$-plumbed Seifert surface is 
called simply \emph{Hopf-plumbed}.

\begin{prop} If $S$ has an annulus $S_0$-presentation 
\emph{(\ref{iterated plumbing})}, then $S$ is a fiber surface 
iff $S_0$ is a fiber surface 
and \emph{(\ref{iterated plumbing})} is a Hopf $S_0$-presentation.  
In particular, an annulus-plumbed surface is a fiber surface iff it 
is Hopf-plumbed. 
\end{prop}
\label{annulus-plumbed fibers iff Hopf-plumbed}
\begin{proof} 
Given that an annulus is a fiber surface iff 
it is a Hopf annulus (\ref{incompressible surfaces}),
the proposition follows immediately from Gabai's theorem \cite{Gabai85} 
that a Murasugi sum of two Seifert surfaces is a fiber surface 
iff both summands are fiber surfaces.
\end{proof}

\section{Basketry and arborescence}
\label{basketry and arborescence}
\subsection{Baskets}\label{basketry}
A \emph{basket $S_0$-presentation} is an annulus $S_0$-presentation 
(\ref{iterated plumbing}) such that each $K_s=O_s$ is an unknot 
and each plumbing arc $\a_s$ is contained in $S_0$.  A Seifert surface 
with a basket $D^2$-presentation is a \emph{basket}.

\begin{prop} A basket is a fiber surface 
iff it is Hopf-plumbed. 
\end{prop}
\begin{proof} Immediate from \ref{annulus-plumbed fibers
iff Hopf-plumbed}. \end{proof}

The failure of (\ref{iterated plumbing}) to distinguish top from bottom	
is, in the case of baskets, alleviated by the following result.

\begin{prop} Any basket is isotopic, by an isotopy
fixing $D^2$, to a basket with the same plumbing arcs 
and plumbands, in which each plumbing is a top-plumbing.
\end{prop}

\begin{proof} For any $S_0$, if $S$ has a basket $S_0$-presentation
with $k$ bottom-plumbed plumbands and $\ell$ top-plumbed plumbands,
then (by \ref{top/bottom}) $S$ has a basket $S_0$-presentation 
in which all bottom-plumbed plumbands precede all top-plumbed 
plumbands, i.e., 
$S=(S_0\BOTPLUMB_{{\vec\a}'}A({\vec O}',{\vec n}'))%
\TOPPLUMB_{{\vec\a}''}A({\vec O}'',{\vec n}'')$, where 
${\vec\a}'=(\a'_1,\dots,\a'_k)$, 
${\vec\a}''=(\a''_{1},\dots,\a''_\ell)$, 
and so on.  

If now $S_0=D^2$ (the basket case), then induction on
$k$ (using \ref{back-to-front} and \ref{top/bottom}) 
completes the proof, and in fact establishes that 
$S \iso_{D^2}\!D^2\TOPPLUMB_{{\vec\a}}A({\vec O},{\vec n})$,
where ${\vec\a}= %
(\a'_k,\a'_{k-1},\dots,\a'_1,\a''_1,\dots,\a''_\ell)$, 
${\vec n}=(n'_k,n'_{k-1},\dots,n'_1,n''_1,\dots,n''_\ell)$.
\end{proof}

\subsection{Arborescent Seifert surfaces}\label{arborescence}
An \emph{arborescent $S_0$-presentation} is an annulus $S_0$-presentation 
(\ref{iterated plumbing}) such that each $K_s=O_s$ is an unknot and 
each plumbing arc $\a_s$ with $s>1$ is a transverse arc 
$\t(\AON{t}\sub S_t)$ for some $t<s$.  A Seifert surface
with an arborescent $D^2$-presentation is simply called
\emph{arborescent}.

\begin{figure}
\begin{center}
\includegraphics{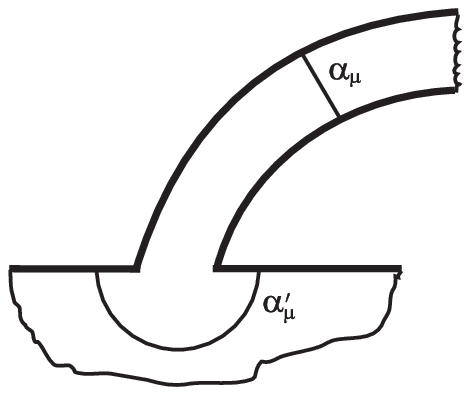}
\end{center}
\caption{Pushing $\t(\AON{t}\sub S'_t)$ off $\AON{t}$.}
\label{Figure 9}
\end{figure}

\begin{prop} \label{arborescent surface is basket}
An arborescent Seifert surface is a basket.
\end{prop}
\begin{proof} More generally and precisely, 
if $S=S_0\PLUMB_{{\vec\a}}A({\vec O},{\vec n})$
is an arborescent $S_0$-presentation with $\mu$ plumbands, 
then there are proper arcs
$\a'_1=\a_1,\a'_2,
\dots, 
\a'_\mu\sub S_0$ such that 
$S':=S_0\PLUMB_{{\vec\a}'}A({\vec O},{\vec n}) \iso_{M_\mu}\!\!\! S$,
where (for suitable transverse arcs and regular neighborhoods)
$M_k:=(S_0\setminus N_{S_0}(\Bd\a_1))\cup%
\bigcup_{s=2}^k N_{\AON{s}}(\t(\AON{s}\sub S))$.

The proof is by induction on $\mu$.  For $\mu\le1$, the assertion 
is trivial.  Let 
$S=(S_0\PLUMB_{{\vec\a}}A({\vec O},{\vec n}))\PLUMB_{\a_\mu} A(O_\mu,n_\mu)$
be an arborescent $S_0$-presentation such that 
$S_{\mu-1}:=S_0\PLUMB_{{\vec\a}}A({\vec O},{\vec n}) \iso_{M_{\mu-1}}\!
S'_{\mu-1}:=S_0\PLUMB_{{\vec\a'}}A({\vec O},{\vec n})$ 
for appropriate proper arcs $\a'_s\sub S_0$.
By assumption, $\a_\mu\sub S_{\mu-1}$ is a transverse arc
$\t(\AON{t}\sub S_t)$ for some $t<s$, so also without loss of generality
$\a_\mu = \t(\AON{t}\sub S'_t)$.
Figure~\ref{Figure 9} illustrates one of
the two possible choices (up to isotopy on $S_0$) 
of $\a'_\mu\sub S_0\sub S'_t\sub S'_{\mu-1}$.
\end{proof}
\begin{remarks}\emph{
(1)~Given $S=D^2\PLUMB_{\vec\a} A(\vec O, \vec n)$,
an arborescent $D^2$-present\-a\-tion of an arborescent 
Seifert surface $S$, it is easy to choose $\vec\a'$ in 
\ref{arborescent surface is basket} so that 
$\As':=\bigcup\a'_s\sub D^2$ is a tree 
such that, if $v\in\Vert(\As')$, then 
$\card\{\eedge\in\Edge(\As'): v\in\Bd\eedge\}\in\{1,4\}$
(for instance, if $\Int D^2$ is identified with the real 
hyperbolic plane, then each $\a_s$ can be taken to be a 
hyperbolic line together with its ideal endpoints on $\Bd D^2$).
There is a partial converse: if $S=D^2\PLUMB_{\vec\a} A(\vec O, \vec n)$ 
is a basket presentation such that $\As:=\bigcup\a_s\sub D^2$ is a
tree and each intrinsic vertex of $\As'$ has valence exactly $4$,
then there is an arborescent surface $\widetilde S$,
having an arborescent presentation with the same plumbands as $S$,
such that $\Bd \widetilde S=\Bd S$ (use 
\ref{back-to-front} repeatedly).
The situation is not merely reminiscent of, but actually
strictly analogous to, the duality between ``resolution trees'' and
``plumbing diagrams'' in the theory of graph manifolds and
resolution of surface singularities (cf.\ \cite{Eisenbud-Neumann},
particularly the figures on p.\ 146).
(2)~Arborescent Seifert surfaces and their boundaries
(there called \emph{special arborescent links}) have also been
studied, from a rather different point of view, in \cite{Sakuma}.
}\end{remarks}

\section{Espaliers, braids, and $\Ts$-bandword 
surfaces}\label{espaliers, braids, etc.}
\subsection{Espaliers}\label{espaliers}
A tree embedded in $\C$ (as a stratified space,
i.e., with all vertices---not just intrinsic vertices---distinguished) 
is \emph{planar}.  A planar tree $\Ts$ is an \emph{espalier} 
if each $\eedge\in\Edge(\Ts)$ is a proper edge in $\LHP$ 
and $\Re|\,\eedge:\eedge\to\R$ is injective.  

\begin{lem} Every planar tree is isotopic to
an espalier.

\end{lem}

Given espaliers $\Ts$ and $\Ts'$, it is useful, 
though abusive, to write $\Ts=\Ts'$ when merely 
$\Vert(\Ts)=\Vert(\Ts')$ and $\Ts$ is isotopic to $\Ts'$ in $\LHP$.
The embedding of an espalier in $\LHP$ is determined 
(up to isotopy in $\LHP$) by the combinatorial structure 
of $\Ts$ (i.e., the underlying abstract simplicial \hbox{$1$-complex}) 
together with the order induced on $\Vert(\Ts)$ by its embedding in $\R$.
In particular, given $n$ real numbers $x_1<\dots<x_n$ and
$n-1$ pairs $\{x_{i(p)},x_{j(p)}\}$ with $1\le i(p)<j(p)\le n$
for $1\le p\le n-1$, the following are equivalent:
(a)~there is an espalier $\Ts$ with 
$\Vert(\Ts)=\{x_1,\dots,x_n\}$, 
$\Edge(\Ts)=\{\eedge_1,\dots,\eedge_{n-1}\}$, and 
$\Bd\eedge_p=\{x_{i(p)},x_{j(p)}\}$;
(b)~for $1\le p<q\le n-1$ the
pairs $(x_{i(p)},x_{j(p)})$ and $(x_{i(q)},x_{j(q)})$ do not
link (i.e., they either touch or unlink).
For $X=\{x_1,\dots,x_n\}\sub\R$, 
$x_1<\dots<x_n$, let
$\Is_X$ (resp., $\Ys_X$) denote any espalier 
$\Ts$ with $\Vert(\Ts)=X$
and 
$\{\Bd\eedge : \eedge\in\Edge(\Ts)\}=\{\{x_p,x_{p+1}\}: 1\le p<n\}$ 
(resp., $\{\{x_1,x_{p}\}: 1< p\le n\}$).
Among the combinatorial types of trees 
with $X$ as the set of $0$-cells, $\Is_X$ and $\Ys_X$ 
represent two extreme types, viz., \emph{linear} (minimal
number of endpoints) and \emph{star-like} (maximal number
of endpoints), respectively; further, among the linear
(resp., star-like) espaliers, $\Is_X$ (resp., $\Ys_X$)
is again extreme, in a sense the reader may formalize
(for $\Ys_X$, see \ref{which T precede Y}).

\subsection{Braid groups; bands}\label{braid groups} 
For $n\in\N$, let $E_n$ denote the \emph{configuration space} 
\begin{displaymath}
\{\{w_1,\dots,w_n\}\sub\C:0\not=\prod\nolimits_{i\ne j}w_i-w_j\}
\end{displaymath}
of unordered $n$-tuples of distinct complex numbers.  
The quotient map of the natural action of $\Ss_n:=\Aut(\numn)$ on 
\begin{displaymath}
\{(w_1,\dots,w_n)\in\C^n:0\not=\prod\nolimits_{i\ne j}w_i-w_j\}
\sub\C^n=\C^{\numn}
\end{displaymath}
induces a topology, smooth structure, and 
orientation on $E_n$, cf.\ \cite{Rudolph83}.
For any basepoint $X\in E_n$, the 
fundamental group $B_X:=\pi_1(E_n;X)$ is \emph{an $n$-string braid group}.
Of course, an $n$-string braid group is isomorphic to the 
\emph{standard} $n$-string braid group $B_n:=B_{\numn}$, 
but it is very convenient to allow more general basepoints.
Denote by $o_X$ the identity of $B_X$.

The map $\{w_1,w_2\}\mapsto(w_1+w_2,(w_1-w_2)^2)$ is 
an orientation-preserving 
diffeomorphism $E_2\to\hbox{$\C\times(\C\setminus\{0\})$}$;
it follows that, for any edge $\eedge\sub\C$, 
the $2$-string braid group $B_{\Bd\eedge}$
is infinite cyclic with a preferred generator, say $\s_\eedge$,
which of course depends only on $\Bd\eedge$.
For $n\ge 2$, if $X\in E_n$ and $X\cap\eedge=\Bd\eedge$, then there is a
natural injection $\iota_{\eedge;X}: B_{\Bd\eedge}\to B_X$,
which depends only on the isotopy class 
of $\eedge$ (rel.\ $\Bd\eedge$) in $\C\setminus(X\setminus\eedge)$.
A \emph{positive $X$-band} 
is any element $\s_{\eedge;X}:=\iota_{\eedge;X}(\s_\eedge)\in B_X$.
(When $X$ is understood, or irrelevant, $\s_{\eedge;X}$ may
be abusively abbreviated to $\s_\eedge$.)
Any two positive $X$-bands are conjugate in $B_X$.
The inverse of a positive $X$-band is a \emph{negative} $X$-band.
Write $|\s_{\eedge;X}^{\pm 1}|:=\s_{\eedge;X}$.

In any group $G$, let $\commute{g}{h}:=ghg^{-1}h^{-1}$ (resp., 
$\YBax{g}{h}:=ghgh^{-1}g^{-1}h^{-1}$) denote 
the \emph{commutator} (resp., the \emph{yangbaxter}) 
of $g, h \in G$.

\begin{lem}\label{relations}
Let $X\in E_n$.  Let $\eedge, \fedge\sub\C$ be 
two edges with $X\cap\eedge=\Bd\eedge, X\cap\fedge =\Bd\fedge$.  
If $\eedge\cap\fedge=\emptyset$ %
{\rm (\emph{resp., $\eedge\cap\fedge=\{x\}\sub \Bd\eedge\cap\Bd\fedge$})}, 
then 
$\commute{\s_{\eedge;X}}{\s_{\fedge;X}}=o_X$ %
{\rm (\emph{resp., $\YBax{\s_{\eedge;X}}{\s_{\fedge;X}}=o_X$})}.
\end{lem}
\begin{proof} Geometrically obvious (trivially so for the commutator, 
slightly less trivially for the yangbaxter).
\end{proof}

\begin{cor}\label{group from tree} Let $\Ts\sub\C$ be a planar tree.
The positive $\Vert(\Ts)$-bands $\s_{\eedge;\Vert(\Ts)}$, 
$\eedge\in\Edge(\Ts)$, generate 
$B_{\Vert(\Ts)}$.  
If $\eedge\cap\fedge=\emptyset$ %
{\rm (\emph{resp., $\eedge\cap\fedge=\{z\}, z\in\Vert(\Ts)$})}, then 
$\commute{\s_{\eedge;\Vert(\Ts)}}{\s_{\fedge;\Vert(\Ts)}}=o_{\Vert(\Ts)}$ %
{\rm (\emph{resp., %
$\YBax{\s_{\eedge;\Vert(\Ts)}}{\s_{\fedge;\Vert(\Ts)}}=o_{\Vert(\Ts)}$})}.
\qed
\end{cor}

Call the $\Vert(\Ts)$-bands $\s_{\eedge;\Vert(\Ts)}$, $\eedge\in\Edge(\Ts)$, 
the \emph{$\Ts$-generators} of $B_{\Vert(\Ts)}$.  

\begin{rem}\label{not the braid group}\emph{
Corollary \ref{group from tree} asserts that the braid group 
$B_{\Vert(\Ts)}$ is a quotient of 
\begin{displaymath}
\operatorname{gp}
\left( 
\s_{\eedge}, \eedge\in\Edge(\Ts)
\left|
\genfrac{}{}{0pt}{}%
{\commute{\s_{\eedge}}{\s_{\fedge}}}
{\YBax{\s_{\eedge}}{\s_{\fedge}}}
\genfrac{}{}{0pt}{}%
{\textrm{ if } \eedge\cap\fedge=\emptyset,}
{\textrm{ if } \eedge\cap\fedge = \{z\}}
\right.
\right);
\end{displaymath}
it does not assert that these groups are identical,
and in fact they are easily seen to be so iff\/ %
$\Ts$ is linear. 
}\end{rem}

\subsection{Embedded band representations 
and $\Ts$-bandwords}\label{embedded band representations}
Let $X\sub\R$ be finite.  If $\eedge\sub\LHP$ is a 
proper edge with $\Bd\eedge=\{p,q\}\sub X$, then 
$\s_{\eedge;X}\in B_X$ depends only on $\{p,q\}$;
write $\s_{p,q;X}:=\s_{\eedge;X}$ and call $\s_{p,q;X}^{\pm 1}$
an \emph{embedded} $X$-band.  
If $X'\sub X$, then there is a well-defined 
injective homomorphism $\iota_{X';X}:B_{X'}\to B_X$ such 
that $\iota_{X';X}(\s_{p,q;X'})=\s_{p,q;X}$ for all $p, q\in X'$ 
with $p\ne q$; this justifies the notation $\s_{p,q}:=\s_{p,q;X}$.  

\begin{prop}\label{more relations} {\rm (1)}~If $\{p,q\}$ 
and $\{r,s\}$ unlink then $\s_{p,q}\s_{r,s}=\s_{r,s}\s_{p,q}$.
{\rm (2)}~If $p<q<r$ then 
$\s_{p,q}\s_{q,r}^{\pm 1}=\s_{p,r}^{\pm 1}\s_{p,q}$.
\end{prop}
\begin{proof} Immediate from \ref{relations}.\end{proof}

An \emph{embedded $X$-band representation} 
is a word $\brep=:(b(1),\dots,b(k))$ 
such that each $b(s)$ is an embedded $X$-band;
for $k=1$ write $\brep=b(1)$ instead of $\brep=(b(1))$.
The \emph{concatenation} of $\brep$
with $\brep'=:(b'(1),\dots,b'(\ell))$ is 
$\brep\concat\brep':=(b(1),\dots,b(k),b'(1),\dots,b'(\ell))$.
The \emph{braid of $\brep$} is 
$\b(\brep):=b(1)\dotsm b(k)\in B_X$.
An embedded $X$-band representation of length $k$ determines 
(and is determined by) a map 
$(i_\brep, j_\brep, \e_\brep): \numk \to X\times X\times\{+,-\}$
with $i_\brep<j_\brep$, such that 
$b(s)=\s_{i_\brep(s),\,j_\brep(s);X}^{\e_\brep(s)1}$.  
Extend $\iota_{X';X}$ to embedded band representations termwise.

Let $\Ts$ be an espalier with $\Vert(\Ts)=X$, so 
the $\Ts$-generators of $B_X$ are embedded $X$-bands.
A \emph{$\Ts$-bandword} is an embedded $X$-band representation
$\brep$ such that each $|b(s)|$ is a $\Ts$-generator.
Let $\brep$ be a $\Ts$-bandword.
Say $\brep$ is \emph{$\eedge$-strict} 
if $\{s: |b(s)|=\s_{\eedge;X}\}\ne\emptyset$, and \emph{strict}
if it is $\eedge$-strict for every $\eedge\in\Edge(\Ts)$.
Say $\brep$ is \emph{$\eedge$-positive}
(resp., \emph{$\eedge$-negative})
if $\brep$ is $\eedge$-strict and 
$\e_\brep|\{s: |b(s)|=\s_{\eedge;X}\}$ is the constant ${+}$
(resp., the constant ${-}$); say 
$\brep$ is \emph{positive} (resp., \emph{homogeneous})
if, for every $\eedge\in\Edge(\Ts)$, 
$\brep$ is $\eedge$-positive (resp., either $\eedge$-positive
or $\eedge$-negative).

\begin{remarks}\emph{
(1)~Embedded $\Is_\numn$-band representations 
are ``embedded band representations in $B_n$'' as defined 
in \cite{Rudolph83}.  The $\Is_\numn$-generator $\s_{i,i+1;\numn}$
is the standard generator $\s_i$ of $B_\numn=B_n$, and 
an $\Is_\numn$-bandword is essentially a \emph{braid word} 
in $B_n$ in the usual sense, \cite[p.\ 70 ff.]{Birman}, 
\cite{Rudolph82}.  The definitions of positive and 
homogeneous $\Ts$-bandwords extend to an arbitrary 
espalier $\Ts$ what is in effect the established usage 
for $\Ts=\Is_\numn$, 
cf.\ \cite{Birman,Stallings,Rudolph83b,Rudolph88}. %
(2)~Of course the set of all $\binom{n}{2}$ positive 
embedded $\Is_\numn$-bands generates $B_n$.
In \cite{Birman-Ko-Lee}, Birman, Ko, and Lee show
that the group with these generators,
and the relations given in \ref{more relations}, 
is in fact $B_n$ (rather than being strictly larger,
cf.\ \ref{not the braid group}); they give 
interesting applications to algorithms for 
the word and conjugacy problems.
}\end{remarks}

\subsection{Braided Seifert surfaces 
and $\Ts$-bandword surfaces}\label{braided Seifert surfaces}
A Seifert surface $S\sub\C\times\R$ is \emph{braided} 
(resp., a \emph{$\Ts$-bandword surface}) 
if $S$ has a handle decomposition (\ref{handle decomposition}) 
with $X, U\sub\R$, satisfying 
\ref{xz-plane separates skeleta}--\ref{saddlepoints}
(resp., \ref{xz-plane separates skeleta}--\ref{lies over espalier}).
\begin{propty}\label{xz-plane separates skeleta}\emph{
$S\cap\UHP\times\R=S^{(0)}$ and 
$S\cap\LHP\times\R=S^{(1)}$.
}\end{propty}

\begin{propty}\label{0-skeleton}\emph{
$\Re(\h0x)=\{x\}$, and 
$\Bd\h0x \cap \{x\}\times\R$ is a single boundary arc of $\Bd\h0x$
on which $\Bd\h0x$ induces the same orientation as that induced
from the orientation of $\{x\}\times\R\sub\C\times\R$.
}\end{propty}

\begin{propty}\label{attaching arcs}\emph{
$u\in\Int(\pr_2(\h1u))\sub\R$, and the closed 
intervals $\pr_2(\h1u)$, $u\in U$, are pairwise disjoint.
}\end{propty}

\begin{propty}\label{saddlepoints}\emph{
$\pr_2|\h1u:\h1u\to\R$ is Morse 
with exactly one (interior) critical point, of index $1$
and critical value $u$, whence $\h1u\cap\LHP\times\{u\}$ is
the union of a transverse arc $\t(\h1u)$ and a core arc $\k(\h1u)$
intersecting transversely at one point (the critical point).
}
\end{propty}

\begin{propty}\label{lies over espalier}\emph{
There is an espalier $\Ts$ with $\Vert(\Ts)=X$
and $\pr_1(\k(\h1u))\in\Edge(\Ts)$ for $u\in U$.
}\end{propty}

Given \ref{xz-plane separates skeleta}--\ref{saddlepoints}
(resp., \ref{xz-plane separates skeleta}--\ref{lies over espalier}),
it is easy to extract an embedded \hbox{$X$-band} representation 
(resp., a $\Ts$-bandword) $\brep_S$:
if $U=:\{u_1,\dots,u_k\}$, $u_1<\dots<u_k$, then set
$\{i_{\brep_S}(s),j_{\brep_S}(s)\}:=\Bd\pr_1(\k(\h1{u_s}))$ and 
$\e_{\brep_S}(s):={\pm}$ if the positive normal vector to $S$ 
at the critical point of $\pr_2|\h1{u_s}$ is a positive 
multiple of $\pm D\pr_2$.  Clearly, $\brep_S=\brep_{S'}$ 
iff $S$ and $S'$ are isotopic through braided Seifert 
surfaces (resp., $\Ts$-bandword surfaces), 
and every band representation (resp., $\Ts$-bandword) 
is $\brep_S$ for some braided Seifert surface (resp., $\Ts$-bandword 
surface) $S$.  Write $S=S(\brep)$ when $\brep=\brep_S$.
Note that a $\Ts$-bandword surface $S=S(\brep)$ determines $\Ts$ 
iff\/ $S$ is connected iff\/ $\brep$ is strict.

\begin{figure}
\begin{center}
\includegraphics{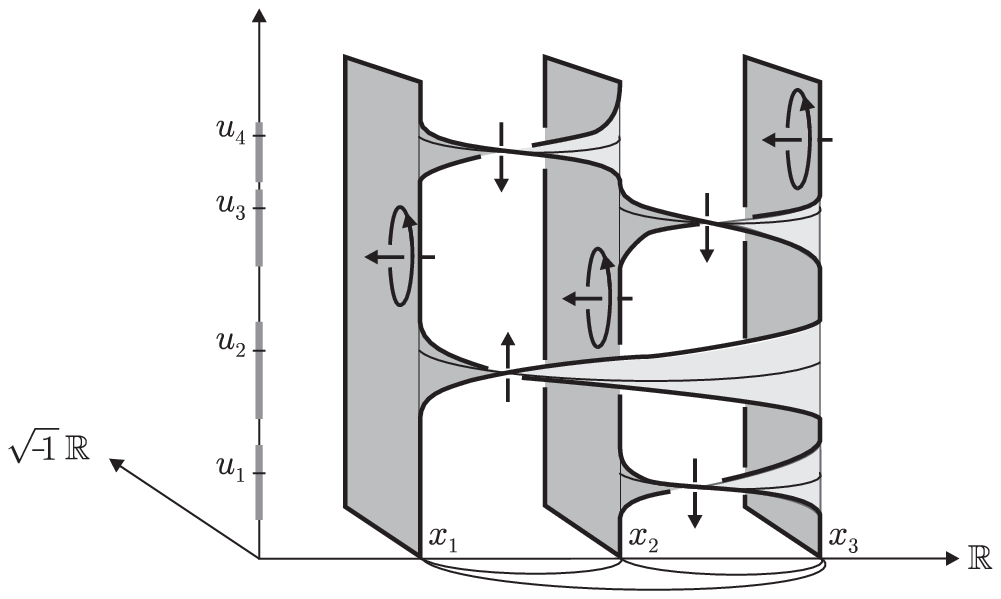}
\end{center}
\caption{A braided surface in $\C\times\R$, 
embellished with normal vectors, 
core arcs, projections of the core arcs into $\LHP\times\{0\}$, 
and projections of the $1$-handles into $\{0\}\times\R$.}
\label{Figure 10}
\end{figure}

Figure~\ref{Figure 10} illustrates a 
braided surface $S(\brep)$ with three \hbox{$0$-handles}
and four \hbox{$1$-handles}, which evidently fails to satisfy
\ref{lies over espalier}.  Deletion of the first
and third, or the second, or the fourth $1$-handle
from $S(\brep)$, produces a $\Ts$-bandword
surface (for various espaliers $\Ts$ 
with $\Vert(\Ts)=\{x_1,x_2,x_3\}$).

For any braided Seifert surface $S(\brep)$, the projection 
$(\Re,\pr_2)(S(\brep) \cap \LHP\times\R)\sub \R\times\R$
becomes a \emph{braid diagram} for $\b(\brep)$ once it is
embellished in the usual way with the correct crossings
at its $k$ doublepoints. 
It is only slightly abusive to conflate $\Bd S(\brep)$ 
with the \emph{closed braid} $\widehat\b(\brep)$.  (In fact,
$\Bd S(\brep)$ and $\widehat\b(\brep)$ have the same link type, 
and it is easy to choose a 
braid axis $A\sub\Int(\UHP)\times\R$ 
and well-positioned $0$-handles $\h0u\sub S(\brep)$ so that 
$\Bd S(\brep)$ actually is a closed braid.)

Call a braided Seifert surface $S$ \emph{standardized} 
if $X=\numn$, $U=\numk$.  
Every isotopy class of braided Seifert surfaces has a 
standardized representative.  It is, however, convenient not to be
limited to standardized braided Seifert surfaces.

\begin{rem}\emph{
Every Seifert surface is isotopic to a 
braided Seifert surface \cite{Rudolph83}.
Not every Seifert surface is isotopic to a $\Ts$-bandword 
surface: e.g., it is easy to see that if $\brep$ is a strict 
$\Ts$-bandword then $\pi_1(S^3\setminus S(\brep);*)$ is free. 
}\end{rem}

\subsection{Elementary moves of braided and
$\Ts$-bandword surfaces}\label{elementary moves}
\begin{figure}
\begin{center}
\includegraphics{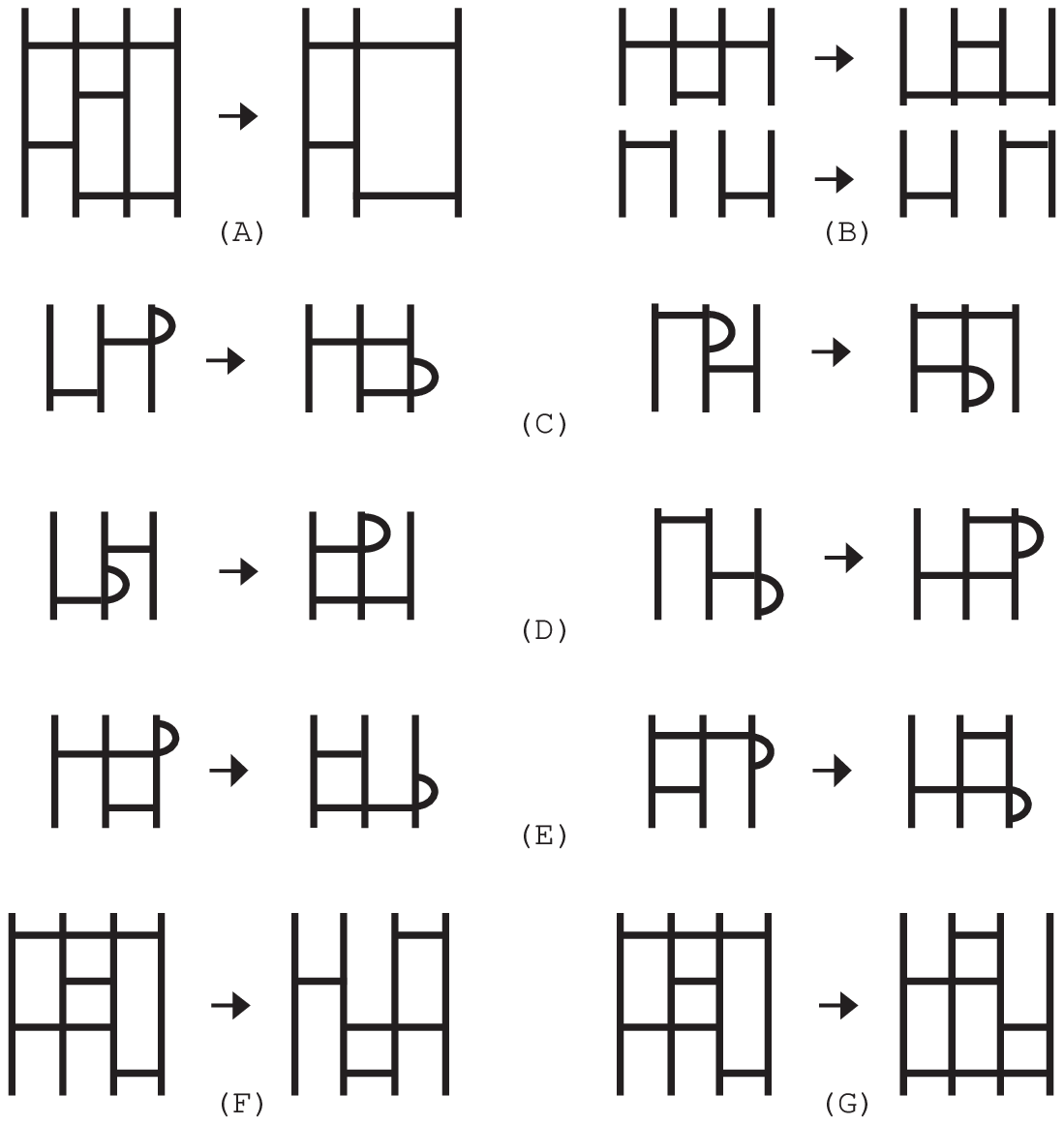}
\end{center}
\caption{(A)~A deflation. (B)~Two slips.
(C)~Two straight slides up.  
(D)~Two straight slides down.
(E)~Two doubly-bent slides up.
(F)~A twirl. (G)~A turn.}
\label{Figure 11}
\end{figure}

Given $S=S(\brep)$ satisfying 
\ref{xz-plane separates skeleta}--\ref{saddlepoints},
if $\Phi_\brep:=(\Re,\pr_2)(G_\brep)\sub \R\times\R$ 
is the projection of 
$G_\brep:=\bigcup_{x\in X}(\h0x\cap\R\times\R)%
\cup \bigcup_{u\in U}\k(\h1u)%
\sub S(\brep)\cap (X\times\R \cup U\times\LHP)$,
then $(\Phi_\brep,\e_\brep)$ is a \emph{charged fence diagram} 
(with \emph{graph} $G_\brep$),
as described in \cite{Rudolph92} and \cite{Rudolph98}
(where $S(\brep)$ was denoted by $S[i,j,\e]$, with $i=i_\brep$
and so on).  Figure~\ref{Figure 11} (adapted 
from \cite{Rudolph98}, Figures 2 and~3)
uses fence diagrams (with the charge on a wire indicated,
where necessary, by a crook at its right end) 
to illustrate several elementary moves, each of which 
replaces a braided Seifert surface $S(\brep)$ by an isotopic 
braided Seifert surface $S(\brep')$ having a different standardization.
These moves can also be described in the language and
notation of embedded $X$-band representations, as follows.

\begin{propty}\label{deflation}\emph{Let $p\in X$, 
$q\in X':=X\setminus\{p\}$.
Let $\brep'_0, \brep'_1$ be embedded $X'$-band 
representations, $\brep':=\brep'_0\concat\brep'_1$.  Say 
$\brep:=\iota_{X';X}(\brep'_0)\concat\s_{p,q;X}^{\pm 1}%
\concat\iota_{X';X}(\brep'_1)$ is obtained from $\brep'$ 
by an \emph{inflation} of sign $\pm$, and that $\brep'$ is 
obtained from $\brep$ by a \emph{deflation}.  Algebraically,
$\b(\brep)$ is obtained from  $\b(\brep')$ by an insignificant
generalization of that ``Markov move'' \cite{Birman} which increases 
braid index by $1$.
}\end{propty}

\begin{propty}\label{slips}\emph{Let $p, q, r, s\in X=:X'$ 
be such that $\{p,q\}$ and $\{r,s\}$ unlink.  
Let $\brep_0, \brep_1$ be embedded 
$X$-band representations.  Say 
$\brep':=\brep_0\concat(\s_{p,q}^{\pm_0 1},\s_{r,s}^{\pm_1 1})\concat \brep_1$ 
is obtained from 
$\brep:=\brep_0\concat(\s_{r,s}^{\pm_1 1},\s_{p,q}^{\pm_0 1})\concat \brep_1$ 
by a \emph{slip}.  Algebraically, 
$\b(\brep)$ is obtained from $\b(\brep')$ by an application 
of \ref{relations}(1).
}\end{propty}

\begin{propty}\label{slides}\emph{Let $p, q, r\in X=:X'$ 
be such that $p<q<r$.  
Let $\brep_0, \brep_1$ be embedded $X$-band representations.
(A)~Say that
$\brep':=\brep_0\concat(\s_{q,r},\s_{p,r}^{\pm 1})\concat \brep_1$ 
(resp., $\brep':=\brep_0\concat(\s_{p,q}^{-1},\s_{p,r}^{\pm 1})\concat %
\brep_1$) is obtained from 
$\brep:=\brep_0\concat(\s_{p,q}^{\pm 1},\s_{q,r})\concat \brep_1$ (resp., 
$\brep:=\brep_0\concat(\s_{q,r}^{\pm 1},\s_{p,r}^{-1})\concat \brep_1$)
by a \emph{straight slide up}.  Call the inverse to a 
straight slide up a \emph{bent slide down}.  (B)~Say 
$\brep':=\brep_0\concat(\s_{p,r}^{\pm 1},\s_{p,q})\concat \brep_1$
(resp., $\brep':=\brep_0\concat(\s_{p,r}^{\pm 1},\s_{q,r}^{-1})\concat %
\brep_1$) is obtained from 
$\brep:=\brep_0\concat(\s_{p,q},\s_{q,r}^{\pm 1})\concat \brep_1$ 
(resp., $\brep:=\brep_0\concat(\s_{q,r}^{-1}, %
\s_{p,q}^{\pm 1})\concat \brep_1$)
by a \emph{straight slide down}.  Call the inverse to a 
straight slide down a \emph{bent slide up}.
(C)~Say 
$\brep':=\brep_0\concat(\s_{p,r},\s_{p,q}^{\pm 1})\concat \brep_1$ 
(resp., $\brep_0\concat(\s_{p,r}^{-1},\s_{q,r}^{\pm 1})\concat \brep_1$) 
is obtained from 
$\brep_0\concat(\s_{q,r}^{\pm 1},\s_{p,r})\concat \brep_1$ 
(resp., $\brep_0\concat(\s_{p,q}^{\pm 1},\s_{p,r}^{-1})\concat \brep_1$) 
by a \emph{doubly-bent slide up}.  Call the inverse to
a doubly-bent slide up a \emph{doubly-bent slide down}.
Algebraically, in all these cases, $\beta(\brep)$ is obtained
from $\beta(\brep')$ by an application of \ref{relations}(2).
}\end{propty}

\begin{propty}\label{twirl}\emph{Let $x'>\max X$, 
$X':=X\setminus\{\min X\}\cup\{x'\}$.
Define $f:X\to X'$ by $f(x)=x$ for $x\in X\cap X'$, $f(\min X)=x'$.
Let $\brep$ be an embedded $X$-band representation.
Define $\brep'$ by setting $\brep'(s)=\s_{f(p),f(q);X'}^{\pm 1}$
when $\brep(s)=\s_{p,q;X}^{\pm 1}$.
Say $\brep'$ is obtained from $\brep$ by a \emph{twirl}.
Algebraically, the standardization of $\beta(\brep')$ (in $B_{\card(X)}$)
is obtained from the standardization of $\beta(\brep)$ by 
conjugation with $\s_{n-1}\s_{n-2}\dotsm\s_1$.
}\end{propty}

\begin{propty}\label{turn}\emph{Let $\brep$ be 
an embedded $X$-band representation,
$X'=X$.
Say $\brep':=(b(k),b(1),\dots,b(k-1))$ is obtained from $\brep$
by a \emph{turn}.
Algebraically, $\beta(\brep')$ is obtained from $\beta(\brep)$ by 
conjugation with the embedded band $b(k)$.
}\end{propty}

\begin{prop}\label{moves preserve isotopy}
If $\brep'$ is obtained from an 
embedded $X$-band representation $\brep$ by 
an inflation, a deflation, slip,
slide (straight, bent, or doubly-bent; 
up or down),
twirl, or turn, then $\brep'$ is 
an embedded $X'$-representation, 
the closed braid $\widehat\beta(\brep')$ is isotopic to 
$\widehat\beta(\brep)$, 
and in fact the braided Seifert surface $S(\brep')$ is 
isotopic to $S(\brep)$.
\end{prop}
\begin{proof} Certainly $\brep'$ is an embedded $X'$-band
representation.  That $\widehat\beta(\brep')$ is isotopic to 
$\widehat\beta(\brep)$ follows, in light of the algebraic 
interpretations of \ref{deflation}--\ref{turn},
from (the easy direction of)
Markov's Theorem \cite{Birman}; of course it also follows from
the final claim, about braided Seifert surfaces, which is 
established by observing that in each case the elementary
move from $\brep'$ to $\brep$ corresponds to an elementary 
move, preserving \ref{xz-plane separates skeleta}--\ref{saddlepoints},
on the ordered handle decomposition (\ref{handle decomposition}) 
of $S(\brep)$.  Specifically: 
an inflation (resp., a deflation) corresponds to adjoining 
(resp., removing) a new $0$-handle and a new $1$-handle
which attaches it to the rest of the surface;
a slip corresponds to transposing the order in
which two adjacent $1$-handles are attached to 
four suitably placed $0$-handles;
the various slides are, precisely, handle slides of one 
$1$-handle over another suitably placed $1$-handle;
and a twirl (resp., a turn)
consists essentially of a cyclic reordering of the index
sets of the $0$-handles (resp., the $1$-handles) 
of (\ref{handle decomposition}).
In turn, each of these moves on (\ref{handle decomposition})
corresponds, in a
standard way, to an isotopy from $S(\brep)$ to $S(\brep')$.
\end{proof}

\begin{rem}\emph{
A more general notion of slide was introduced 
in \cite{Rudolph83}, in the context of Seifert ribbons
and their corresponding (not necessarily embedded)
band representations.  The straight, bent, and doubly-bent 
slides defined above are those which preserve embeddedness 
and so are suitable in the present context of Seifert surfaces.
}\end{rem}

\begin{hist}\label{history of fences}\emph{
It should have been 
mentioned in \cite{Rudolph92} or \cite{Rudolph98} that 
the first published appearance of fences 
(and certainly the first time I saw them)
may well have been c.\ 1959, in one of Martin Gardner's ``Mathematical 
Games'' columns in \emph{Scientific American}, reprinted as Chapter~2, 
``Group Theory and Braids'', of \cite{Gardner}.  Gardner
refers to them there as ``vertical lines and shuttles'', and 
includes, as an application, a drinking game for computer programmers.
}\end{hist}

\subsection{Moves of espaliers; more isotopies of 
$\Ts$-bandword surfaces}\label{moves of espaliers} 
Let $\Ts$ be an espalier.  Several elementary moves,
with obvious similarities to the
identically named moves on braided Seifert surfaces
(or bandwords) described in \ref{elementary moves},
may be applied to $\Ts$ to produce another espalier $\Ts'$.

\begin{propty}\label{espalier deflation}\emph{
Let $\eedge\in\Edge(\Ts)$ be such that 
$q\in\Bd\eedge\cap\Endpt(\Ts)$.
Let $\Vert(\Ts'):=\Vert(\Ts)\setminus\{q\}$, 
$\Edge(\Ts'):=\Edge(\Ts)\setminus\{\eedge\}$.
Say $\Ts'$ is obtained from $\Ts$ by a \emph{deflation}
and that $\Ts$ is obtained from $\Ts'$ 
by an \emph{inflation}.
}\end{propty}

\begin{propty}\label{esplaier slides}\emph{
Let $\eedge_0\in\Edge(\Ts)$, $\Bd\eedge_0=:\{p,q\}$,
$\{\eedge_0,\eedge_1,\dots,\eedge_r\}:=\{\eedge\in\Edge(\Ts) %
: q\in\Bd\eedge\}$.  
There exist proper edges $\eedge'_1,\dots,\eedge'_r\sub\LHP$
with $\Bd\eedge'_s=\Bd\eedge_s\setminus\{q\}\cup\{p\}$, $s=1,\dots,r$,
such that, if $\Vert(\Ts'):=\Vert(\Ts)$, 
$\Edge(\Ts'):=\Edge(\Ts)\setminus\{\eedge_1,\dots,\eedge_r\}%
\cup\{\eedge'_1,\dots,\eedge'_r\}$, then $\Ts'$ is an espalier.
If $p<q$ (resp., $q<p$), then say $\Ts'$ is obtained 
from $\Ts$ by a \emph{slide left} (resp., a \emph{slide right})
\emph{along $\eedge_0$},
and---of course---that $\Ts$ is obtained from $\Ts'$ by a slide right 
(resp., a slide left) along $\eedge_0$.
}\end{propty}

\begin{propty}\label{espalier twirl}\emph{
Let $q>\max \Vert(\Ts)$.  Let 
$\{\eedge_1,\dots,\eedge_r\}:=\{\eedge\in\Edge(\Ts) : %
\min\Vert(\Ts)\in\Bd\eedge\}$.
There exist proper edges $\eedge'_1,\dots,\eedge'_r\sub\LHP$
with $\Bd\eedge'_s=\Bd\eedge_s\setminus\{\min\Vert(\Ts)\}\cup\{q\}$, 
$s=1,\dots,r$, such that, if $\Vert(\Ts'):=\Vert(\Ts)$, 
$\Edge(\Ts'):=\Edge(\Ts)\setminus\{\eedge_1,\dots,\eedge_r\}%
\cup\{\eedge'_1,\dots,\eedge'_r\}$, then $\Ts'$ is an espalier.
Say $\Ts'$ is obtained from $\Ts$ by a \emph{twirl}.
}\end{propty}

\begin{lem}\label{change of espalier}
If $\brep'$ is obtained 
from a $\Ts$-bandword $\brep$ by an inflation (resp., 
a deflation; a slip; a twirl; a turn),
then $\brep'$ is a $\Ts'$-bandword, where $\Ts'$ is obtained from
$\Ts$ by an inflation (resp., a deflation; doing nothing; a twirl; 
doing nothing).
\qed
\end{lem}

\begin{figure}
\begin{center}
\includegraphics{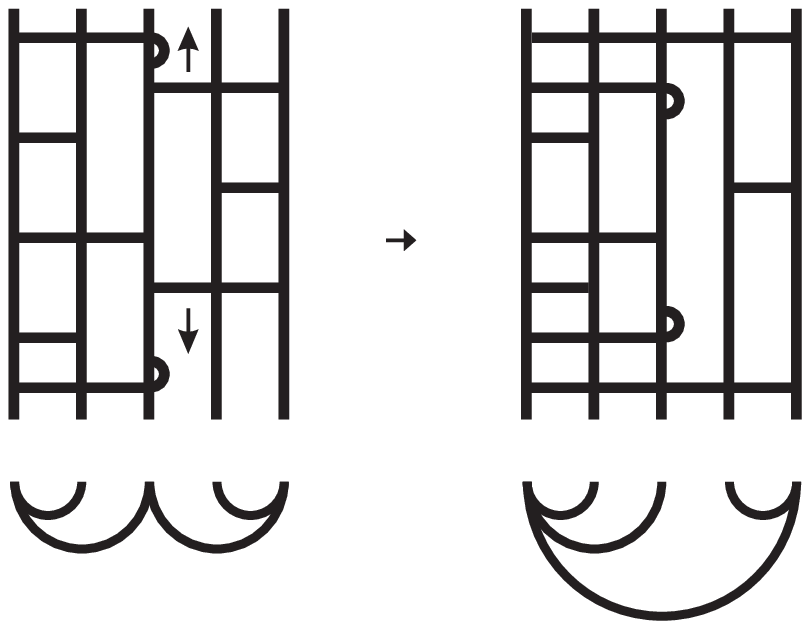}
\end{center}
\caption{Above, partially charged fence diagrams for a $\Ts$-bandword 
surface $S(\brep)$ and a $\Ts'$-bandword surface $S(\brep')$
isotopic to $S(\brep)$ by a slide up, a slip, 
and a slide down; below, $\Ts$ and $\Ts'$.}
\label{Figure 12}
\end{figure}

In contradistinction to \ref{change of espalier}, 
if $\brep'$ is obtained from $\brep$ by 
applying a slide (straight, bent, or doubly-bent; up or down), then
(with trivial exceptions) the embedded $\Vert(\Ts)$-band representation 
$\brep'$ is not a $\Ts'$-bandword for any espalier $\Ts'$, 
in particular not for any $\Ts'$ produced from $\Ts$ by a slide.
It is nevertheless the case that for certain $\Ts$-bandwords 
$\brep$, if an appropriate sequence of two or more 
slides, with slips interspersed as necessary, is applied to $\brep$,
then the result is a $\Ts'$-bandword $\brep'$, where 
$\Ts'$ is produced from $\Ts$ by a slide.
Figure~\ref{Figure 12} illustrates such a multiple 
slip-slide of bandword surfaces covering a slide of espaliers;
on each side of the figure, the upper part is a fence diagram 
$(\Re,\pr_2)(G)$ (with slight distortions, as in 
Figure~\ref{Figure 11}, 
to indicate charges)) derived from the graph $G$ on the bandword
surface, while the lower part is the espalier $\pr_1(G)$.
\section{Proof of the Main Theorem}\label{characterization}

\subsection{Characterization of Hopf-plumbed 
baskets}\label{homogeneous=Hopf-plumbed}

Write $\Ts\prec \Ts'$ 
in case $\Ts$ and $\Ts'$ are espaliers with $\Vert(\Ts)=\Vert(\Ts')$
and either 
\begin{displaymath}
\card(\Endpt(\Ts))\!<\!\card(\Endpt(\Ts'))
\end{displaymath}
or	
\begin{displaymath}
\card(\Endpt(\Ts))=\card(\Endpt(\Ts')) {\textrm{ and }}
\sum_{v\in\Endpt(\Ts)} v < \sum_{v\in\Endpt(\Ts')} v .
\end{displaymath}
Let $\Ts$ be an espalier, $\brep$ a $\Ts$-bandword.
Let $\Ts$ be an espalier, $\brep$ a $\Ts$-bandword.
Recall the definition of $\Ys_X$ from \ref{espaliers}.
The following lemma formalizes the assertion there 
that $\Ys_X$ is the most extreme star-like espalier 
with $X$ as the set of $0$-cells.

\begin{lem}\label{which T precede Y}
$\Ts\prec\Ys_{\Vert(\Ts)}$ iff 
$\Ts\ne \Ys_{\Vert(\Ts)}$ iff $\max\Bd\eedge\notin\Endpt(\Ts)$
for some $\eedge\in\Edge(\Ts)$.
\qed
\end{lem}

\begin{lem}\label{move T towards Y}
If $\Ts\ne\Ys_{\Vert(\Ts)}$ 
and $\brep$ is homogeneous, 
then there exists $\Ts'$ with $\Ts\prec\Ts'$ 
and a homogeneous $\Ts'$-bandword $\brep'$ 
with $S(\brep')\iso S(\brep)$.
\end{lem}

\begin{figure}
\begin{center}
\includegraphics{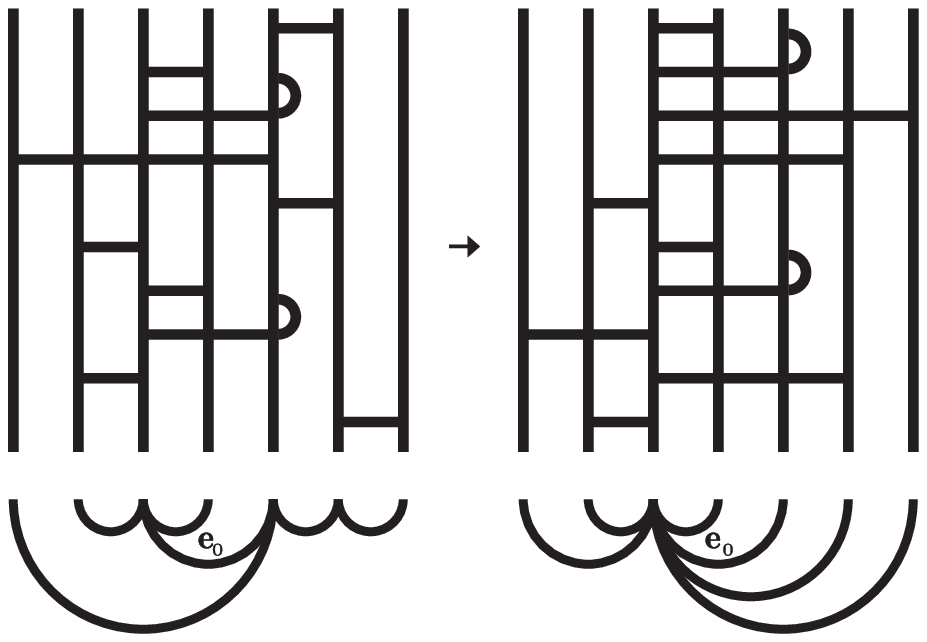}
\end{center}
\caption{Above, partially charged fence diagrams for an 
$\eedge_0$-positive $\Ts$-bandword surface $S(\brep)$ and 
an $\eedge_0$-positive $\Ts'$-bandword surface $S(\brep')$ 
isotopic to $S(\brep)$ by an inverse turn, six slips, and 
four slides down; below, $\Ts$ and $\Ts'$.
}
\label{Figure 13}
\end{figure}

\begin{proof} By \ref{which T precede Y}, there 
exists $\eedge_0\in\Edge(\Ts)$ 
such that $v:=\max\Bd\eedge_0\notin\Endpt(\Ts)$ 
and for all $\eedge\in\Edge(\Ts)\setminus\{\eedge_0\}$, 
if $\max\Bd\eedge=v$, then $\min\Bd\eedge<\min\Bd\eedge_0=:u$.  
For $x\in\{u,v\}$, let $(\Ts-\eedge_0)_x$ denote the component of   
$\Ts\setminus\Int\eedge_0$ which contains $x$.
Clearly $(\Ts-{\eedge_0})_u$ and $(\Ts-{\eedge_0})_v$ 
are disjoint espaliers with 
$\Edge((\Ts-{\eedge_0})_u)\cup\Edge((\Ts-{\eedge_0})_v)=%
\Edge(\Ts)\setminus\{\eedge_0\}$, and 
$\commute{\s_\fedge}{\s_\gedge}=o_{\Vert(\Ts)}$
for all $\fedge\in\Edge((\Ts-{\eedge_0})_u)$ 
and $\gedge\in\Edge((\Ts-{\eedge_0})_v)$.  

In case $\brep$ is $\s_{\eedge_0}$-positive,
an appropriate sequence of turns (\ref{turn}) 
and slips (\ref{slips})
converts $\brep$ to a homogeneous $\Ts$-bandword
$\brep''=\drep_1\concat\drep_2\concat\dotsm\concat\drep_\ell$,
where $\ell\ge1$ and for $1\le s\le\ell$,
$\drep_s=\s_{\eedge_0}\concat\arep_s\concat\crep_s$
with $\arep_s$ a $(\Ts-{\eedge_0})_v$-bandword and 
$\crep_s$ a $(\Ts-{\eedge_0})_u$-bandword; 
by \ref{moves preserve isotopy}, $S(\brep'')\iso S(\brep)$.
By the choice of $\eedge_0$, if 
$\fedge\in\Edge((\Ts-{\eedge})_v)$
and $v\in\Bd\fedge$ then either 
$\min\Bd\fedge=v$ or $\min\Bd\fedge<u$, 
so a further sequence of slips and (straight and bent) 
slides down (\ref{slides}) converts $\brep''$ to a homogeneous 
$\Ts'$-bandword $\brep'$, with $S(\brep')\iso S(\brep)$ 
(again by \ref{moves preserve isotopy}), 
where $\Ts'$ is produced from $\Ts$ by a 
slide right along $\eedge_0$.  
(Figure~\ref{Figure 13} illustrates such 
a move from $\brep$ to $\brep'$, in the style 
of Figure~\ref{Figure 12}.)
Now, if $u\in\Endpt(\Ts)$ 
(i.e., if $(\Ts-\eedge)_u=\{u\}$), then 
$\Endpt(\Ts')=\Endpt(\Ts)\setminus\{u\}\cup\{v\}$,
whereas if $u\notin\Endpt(\Ts)$ then 
$\Endpt(\Ts')=\Endpt(\Ts)\cup\{v\}$; so
in any case $\Ts\prec \Ts'$.

In case $\brep$ is $\s_{\eedge_0}$-negative, 
the procedure is entirely similar, with 
$\drep_s=\arep_s\concat\crep_s\concat\s_{\eedge_0}^{-1}$,
and $\arep_s, \crep_s$ as before; 
slides up take the place of slides down.  
The final products $\brep'$ and $\Ts'$ are as in the first case;
again $S(\brep')\iso S(\brep)$ and $\Ts\prec\Ts'$.
\end{proof}

\begin{cor}\label{from any T to Y} 
If $\brep$ is homogeneous, then there is a homogeneous 
$\Ys_{\Vert(\Ts)}$-bandword $\brep'$ with $S(\brep')\iso S(\brep)$.
\qed\end{cor}

For any $\Ts$-bandword $\arep$, let $c(\arep):=
\sum_{\eedge\in\Edge(\Ts)}\left|{\card(\{s:|a(s)|=\s_\eedge\})-2}\right|$.

\begin{lem}\label{simplify Y-bandwords}
If $\Ts=\Ys_X$ and $\brep$ is homogeneous, 
then for some $X^*$ there is a homogeneous $\Ys_{X^*}$-bandword $\brep^*$ 
such that $S(\brep^*)\iso S(\brep)$ and $c(\brep^*)=0$.
\end{lem}
\begin{proof} 
Since $\brep$ is strict, if $c(\brep)>0$ then either 
$\card(\{s:|b(s)|=\s_{\eedge_0}\})=1$ for some $\eedge_0\in\Edge(\Ys_X)$,
or $\card(\{s:|b(s)|=\s_{\eedge_0}\})\ge 3$ for some 
$\eedge_0\in\Edge(\Ys_X)$.  In either case,
let $v:=\max\Bd\eedge_0$, $u:=\min\Bd\eedge_0=\min X$.

In the first case, there is a deflation (\ref{deflation}) of 
$\brep$ to a homogeneous $\Ys_{X'}$-bandword 
$\brep'$ with $c(\brep')=c(\brep)-1$, where 
$X':=X\setminus\{v\}$; by \ref{moves preserve isotopy}, 
$S(\brep')\iso S(\brep)$.
In the second case, choose $w\in {]v,\min(X\cap{]v,\infty[})[}$
and set $X':=X\cup\{w\}$.  
If $\brep$ is $\eedge_0$-positive, then 
$\brep=\brep_0\concat\s_{\eedge_0}\concat%
\brep_1\concat\s_{\eedge_0}\concat\brep_2\concat%
\s_{\eedge_0}\concat\brep_3$
where $\s_{\eedge_0}=\s_{u,v}$ does not appear in 
$\brep_1$ or $\brep_2$,
and the sequence of moves 
\begin{align*}
\brep & \to\brep_0\concat\s_{u,v}\concat\brep_1\concat%
        \s_{u,v}\concat\s_{v,w}%
        \concat\brep_2\concat\s_{u,v}\concat\brep_3
                    && \text{by (\ref{deflation})} \\
      & \to\brep_0\concat\s_{u,v}\concat\brep_1\concat%
        \s_{v,w}\concat\s_{u,w}%
        \concat\brep_2\concat\s_{u,v}\concat\brep_3
                    && \text{by (\ref{slides})} \\
      & \to\brep_0\concat%
        \s_{u,v}\concat\s_{v,w}%
        \concat\brep_1\concat\s_{u,w}%
        \concat\brep_2\concat\s_{u,v}\concat\brep_3
                    && \text{by (\ref{slips})} \\
      & \to\brep_0\concat%
        \s_{u,w}\concat\s_{u,v}%
        \concat\brep_1\concat\s_{u,w}%
        \concat\brep_2\concat\s_{u,v}\concat\brep_3:=\brep'
                    && \text{by (\ref{slides})}
\end{align*}
converts $\brep$ to a homogeneous $\Ys_{X'}$-bandword $\brep'$
with $c(\brep')=c(\brep)-1$, and $S(\brep') \iso S(\brep)$
by \ref{moves preserve isotopy};
if $\brep$ is $\eedge_0$-negative,
a similar sequence of moves (an inflation of sign $-$, 
a straight slide down, slips, and a straight slide up)
has the same effect.
Induction on $c(\brep)$ completes the proof.
\end{proof}

\begin{lem}\label{Y-homogeneous is Hopf-plumbed basket}
If $\Ts=\Ys_X$, $\brep$ is homogeneous,
and $c(\brep)=0$, then $S(\brep)$ is a Hopf-plumbed basket with
$\card(X)-1$ plumbands. 
\end{lem}
\begin{proof} Without loss of generality, 
take $S(\brep)$ to be standardized, so that
$X=\numn$ and $U=\numk$.  
For $t=2,\dots,n$, write $u(t):=\min j_\brep^{-1}(t)$,
$v(t):=\max j_\brep^{-1}(t)$, and let
$\a_t$ be a proper arc in $\h01$ joining
$\k(\h1{u(t)})\cap\Bd\h01$ to 
$\k(\h1{v(t)})\cap\Bd\h01$.
Evidently 
$\h1{u(t)}\cup\h0t\cup\h1{v(t)}\cup\Nb{{\h01}}{\a_t}$
is a Hopf annulus $A(O_t,\phi(t))$, where $\phi(t)=-1$ 
(resp., $1$) if $\brep$ is $\s_{1,t}$-positive 
(resp., $\s_{1,t}$-negative), and 
$S(\brep)$ is a Hopf-plumbed $\h01$-basket with
plumbing arcs $\a_t$ and plumbands $A(O_t,\phi(t))$.  
\end{proof}

\begin{thm}\label{characterization theorem} 
\emph{(A)}~If $\Ts$ is an espalier
and $\brep$ is a homogeneous $\Ts$-bandword, 
then $S(\brep)$ is a Hopf-plumbed basket.  
\emph{(B)}~If $S$ is a Hopf-plumbed basket, then $S$ is
isotopic to $S(\brep)$ for some espalier $\Ts$ and
homogeneous $\Ts$-bandword $\brep$.
\end{thm}
\begin{proof} (A) follows directly from 
\ref{move T towards Y}-\ref{Y-homogeneous is Hopf-plumbed basket}.
The construction in the proof of 
\ref{Y-homogeneous is Hopf-plumbed basket}, run in reverse, 
proves a more precise formulation of (B): 
if $S$ is a Hopf-plumbed basket with $n-1$ plumbands, 
then $S$ is isotopic to $S(\brep)$, where $\brep$ is
a homogeneous $\Ys_\numn$-bandword with $c(\brep)=0$.
\end{proof}

\subsection{The rest of the Main Theorem}\label{main theorem}
Again let $\Ts$ be an espalier, $\brep$ a $\Ts$-bandword,
and consider the handle decomposition 
(\ref{handle decomposition}) of $S(\brep)$,
satisfying \ref{xz-plane separates skeleta}--\ref{lies over espalier},
with $X=\Vert(\Ts)$.

\begin{lem}\label{not strict}
If $\brep$ is not strict, then 
$S(\brep)$ is not connected. 
\end{lem}
\begin{proof} If $\eedge\in\Edge(\Ts)$ and 
$\brep$ is not $\eedge$-strict, then the $0$-handles 
$\h0{\min\Bd\eedge}$ and $\h0{\max\Bd\eedge}$ must 
lie in different components of $S(\brep)$.
\end{proof}

\begin{figure}
\begin{center}
\includegraphics{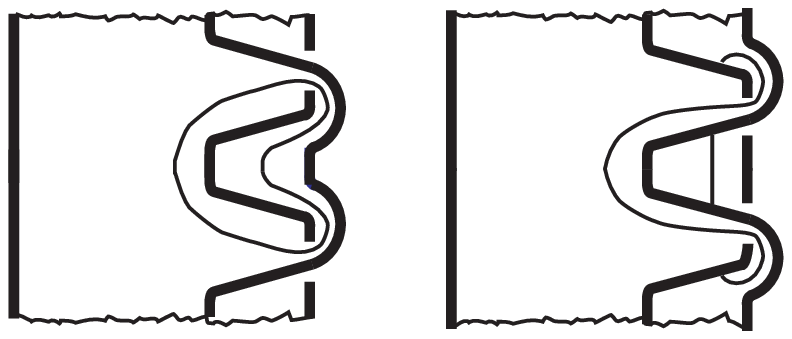}
\end{center}
\caption{The boundaries of top-compression disks for 
$S(\dots,\s{\vrule height 8pt depth 3pt width0pt}_1,%
\s{\vrule height 8pt depth 3pt width0pt}_1^{-1},\dots)$ 
and for $-S(\dots,\s{\vrule height 8pt depth 3pt width0pt}_1^{-1},%
\s{\vrule height 8pt depth 3pt width0pt}_1,%
\dots)$.}
\label{Figure 14}
\end{figure}

\begin{lem}\label{compression (two strings)}
If $X=\mathbf 2$, so $\Ts=\Ys_{\mathbf 2}$,
and $\brep$ is strict but not homogeneous, then both $S(\brep)$ and
$-S(\brep)$ are top-compressible. \end{lem}
\begin{proof} Let $\brep=(\sigma_1^{\e(1)1},\dots,\sigma_1^{\e(k)1})$,
$\e(s)\in\{{+},{-}\}$, be strict but not homogeneous, so that  
$k\ge2$ and there exists $s$ with $\e(s)=-\e(s+1)$.
If $\e(s)={+}$ and $\e(s+1)={-}$
(resp., $\e(s)={-}$ and $\e(s+1)={+}$),
then the union of $\k(\h1s)$, $-\k(\h1{s+1})$, and 
suitable proper arcs on $\h01$ and $\h02$
is the boundary of a top-compression disk for $S(\brep)$ 
(resp., $-S(\brep)$); Figure~\ref{Figure 14} illustrates these
two cases.  Up to turns (\ref{turn}) (which do not change the 
isotopy type of $S(\brep)$, \ref{moves preserve isotopy}), 
both cases occur, and so both $S(\brep)$ and $-S(\brep)$ are
top-compressible.\end{proof}

Let $\Xend:=\Endpt(\Ts)$.  Call $\eedge\in\Edge(\Ts)$ a 
\emph{terminal edge} of $\Ts$ if $\Bd\eedge\cap\Xend\ne\emptyset$.  

\begin{lem}\label{compression (n strings)}
If $\Ts$ is arbitrary, 
and $\brep$ is strict but not homogeneous, 
then $S(\brep)$ is compressible. \end{lem}
\begin{proof} Let $n=\card(X)$.  The case $n=1$ is vacuous.
The case $n=2$ follows from \ref{compression (two strings)}.  
Suppose $n>2$.  
There exists a terminal edge $\eedge$ of $\Ts$ with 
${]}\min\Bd\eedge,\max\Bd\eedge{[}\,\cap X=\emptyset$.
Let $X':=\Bd\eedge$, $X'':=X\setminus(X'\setminus\Xend(\Ts))$,
$U':=\{u\in U: \pr_1(\k(\h1u))=\eedge\}$, $U'':=U\setminus U'$,
$S'=\bigcup_{x\in X'}\h0x\cup \bigcup_{u\in U'}\h1u$,
$S''=\bigcup_{x\in X''}\h0x\cup \bigcup_{u\in U''}\h1u$.
Clearly $S'\cup S''=S(\brep)$ and $S'\cap S''=\h0q$, where 
$\Bd\eedge=:\{p,q\}$, $p\in\Xend(\Ts)$, $q\notin\Xend(\Ts)$;
further, if $q=\max\Bd\eedge$ (resp., $q=\min\Bd\eedge$), 
then (for appropriate choices of collars)
$S'\setminus\h0q\sub\Int\top(S'')$ and 
$S''\setminus\h0q\sub\Int\bot(S')$
(resp., $S'\setminus\h0q\sub\Int\bot(S'')$ and 
$S''\setminus\h0q\sub\Int\top(S')$).
Each of $S'$, $S''$ is evidently a braided Seifert surface 
in its own right: more precisely, 
$S':=S(\iota_{X';X}(\crep))$ 
and $S'':=S(\iota_{X'';X}(\drep))$,
where $\crep$ is an $\eedge$-bandword and
$\drep$ is a $(\Ts-\eedge)_q$-bandword.
If $\brep$ is strict but not homogeneous, 
then both $\crep$ and $\drep$ are strict,
and at least one of $\crep$ and $\drep$ is not homogeneous.

If $\crep$ is not homogeneous, and 
$q=\max\Bd\eedge$ (resp., $q=\min\Bd\eedge$), 
then (by \ref{not strict}) 
there exists a top-compression 
disk for $S'$ (resp., $-S'$), which can be taken 
to be disjoint from $S''$, and therefore to be
a top-compression disk for $S(\brep)$ 
(resp., $-S(\brep)$).  

If $\crep$ is homogeneous
and $\drep$ is inhomogeneous, then $S'$ is a 
Hopf-plumbed basket (by \ref{Y-homogeneous is 
Hopf-plumbed basket}) and $S''$ is
compressible (by induction on $n$).  In fact, 
the proof of \ref{simplify Y-bandwords}
is easily adapted to show that 
$S(\brep)$ has a Hopf basket $S''$-presentation.
By \ref{indifference}, if there exists a 
top-compression disk for $S''$ (resp., $-S''$), 
then there exists such a disk which 
is disjoint from $S'\setminus S''$,
and is therefore also a top-compression disk 
for $S(\brep)$ (resp., $-S(\brep)$).

In any case, if $\brep$ is strict but not homogeneous
then $S(\brep)$ is compressible.\end{proof}

\begin{thm}\label{Main Theorem} The following are equivalent:
\emph{(A)}~$\brep$ is homogeneous;
\emph{(B)} $S(\brep)$ is a Hopf-plumbed basket;
\emph{(C)}~$S(\brep)$ is a fiber surface;
\emph{(D)}~$S(\brep)$ is incompressible and connected.  
\end{thm}
\begin{proof} That (A) implies (B) is contained in 
Theorem \ref{characterization}.
It is a standard fact (\cite{Stallings,Harer,Morton,Gabai83}; 
\ref{annulus-plumbed fibers iff Hopf-plumbed}) 
that any Hopf-plumbed surface is a fiber surface,
so in particular (B) implies (C).  
It is likewise standard (\cite{Stallings,Gabai83b}; 
\ref{incompressible surfaces})
that any fiber surface is incompressible and connected,
so in particular (C) implies (D).
Finally, by \ref{not strict} and 
\ref{compression (n strings)}, not-(A) implies not-(D).
\end{proof}

\providecommand{\bysame}{\leavevmode\hbox to3em{\hrulefill}\thinspace}

\end{document}